\documentclass[12pt,a4paper]{article}

\usepackage{latexsym}
\usepackage{amssymb,amsmath}
\usepackage{graphicx}
\usepackage{ae}
\usepackage{colortbl}
\usepackage{tikz,xtab}
\usetikzlibrary{decorations.pathreplacing}

 \def\f{\frac} \def\eot{\hfill{$\Box$}}
\newcommand{\lettnum}{\def\labelenumi{(\alph{enumi})}} 
\newtheorem{thm}{Theorem}[section]
\newtheorem{lemma}[thm]{Lemma}
\newtheorem{defn}[thm]{Definition}

\newtheorem{cor}[thm]{Corollary}

\tikzstyle{every node}=[circle, draw, fill=black!50,
                        inner sep=0pt, minimum width=4pt]

\begin{document} 
\title{Supercards, Sunshines and Caterpillar Graphs}
\author{Paul Brown\footnote{paulb@dcs.bbk.ac.uk} $\,$ and 
Trevor Fenner\footnote {trevor@dcs.bbk.ac.uk}
\\Department of Computer Science and Information Systems
\\Birkbeck College, University of London, United Kingdom}
\maketitle
\vspace{-10mm}
\begin{abstract}
{\rm The vertex-deleted subgraph $G-v$, obtained from the graph $G$ by deleting the vertex $v$ and all edges incident to $v$, is called a {\em card} of $G$. The {\em deck} of $G$ is the multiset of its unlabelled cards. The {\em number of common cards} $b(G,\,H)$ of $G$ and $H$ is the cardinality of the multiset intersection of the decks of $G$ and $H$. A {\em supercard} $G^+$ of $G$ and $H$ is a graph whose deck contains at least one card isomorphic to $G$ and at least one card isomorphic to $H$. We show how maximum sets of common cards of $G$ and $H$ correspond to certain sets of permutations of the vertices of a supercard, which we call {\em maximum saturating sets}. We apply the theory of supercards and maximum saturating sets to the case when $G$ is a sunshine graph and $H$ is a caterpillar graph. We show that, for large enough $n$, there exists some maximum saturating set that contains at least $b(G,\, H)-2$ automorphisms of $G^+$, and that this subset is always isomorphic to either a cyclic or dihedral group. We prove that $b(G,\,H) \le \f{2(n+1)}{5}$ for large enough $n$, and that there exists a unique family of pairs of graphs that attain this bound. We further show that, in this case, the corresponding maximum saturating set is isomorphic to the dihedral group.}
\end{abstract}
\maketitle
\section{\normalsize Introduction}\label{introduction}\vspace{-5mm}
In this paper all graphs are finite, undirected and contain no loops or multiple edges. Any graph-theoretic terminology and notation not explicitly explained below can be found in Bondy and Murty's text \cite{b&m1}. For more information on the action of a permutation group on the vertices of a graph, we refer the reader to the book by Lauri and Scapellato \cite{lau1}.

Let $G$ be a graph of order $n$ and let $u, v \in V(G)$. We denote the group of all permutations of $V(G)$ by $S_{V(G)}$ and the identity permutation of $S_{V(G)}$ by $1_{V(G)}$. A {\em transposition} of $V(G)$ is a permutation that swaps two vertices in $V(G)$ and leaves the rest unchanged.

The {\em neighbourhood} of $v$ in $G$ is the set $N_G(v)$ consisting of all vertices of $G$ adjacent to $v$. The cardinality of this set is the {\em degree} of $v$ in $G$, i.e., $d_G(v)=|N_G(v)|$. A {\em leaf} of $G$ is a vertex of degree $1$, and an {\em isolated vertex} of $G$ is a vertex of degree $0$. We denote the number of vertices of degree $k$ in $G$ by $d_k(G)$, so $\sum_i d_i(G)=n$. A {\em non-trivial} component of a graph is one of order at least two.

Suppose that $H$ is another graph and that $\gamma$ is a bijection from $V(G)$ to $V(H)$. For any $Z \subseteq V(G)$, we write the image of $Z$ under $\gamma$ as $\gamma(Z)$. When $\gamma$ is, moreover, an isomorphism from $G$ to $H$, i.e., $xy$ is an edge of $G$ if and only if $\gamma(x)\gamma(y)$ is an edge of $H$, we write $\gamma(G)=H$. We write $G \cong H$ to indicate that $G$ and $H$ are isomorphic. The group of all automorphisms of $G$, i.e., isomorphisms from $G$ to itself is denoted by $\mathrm{Aut}(G)$. We note that any transposition of $V(G)$ that swaps a pair of leaves adjacent to the same vertex is in $\mathrm{Aut}(G)$.

Now let $Z \subseteq V(G)$. The $Z$-deleted subgraph $G-Z$ is obtained from $G$ by deleting all the vertices of $Z$ together with all edges of $G$ incident to a vertex in $Z$. So $d_{G-Z}(v)=d_G(v)-|N_G(v) \cap Z|$, for all $v \in V(G-Z)$. When $Z=\left\{v\right\}$ or $Z=\left\{u,\, v\right\}$, we write $G-Z$ as $G-v$ or $G-u-v$, respectively. The {\em vertex-deleted subgraph} $G-v$ is also known as a {\em card} of $G$, and the multi-set of all $n$ unlabelled cards of $G$ is called the {\em deck} of $G$, which we denote by $\mathcal D(G)$. If $Z$ is the set of leaves and isolated vertices of $G$ then $G-Z$ is called the {\em skeleton} of $G$, denoted by $\mathrm {skel}(G)$.

Clearly, if $G \cong H$ then $\mathcal D(G)=\mathcal D(H)$. The {\em Reconstruction Conjecture}, first proposed by Kelly and Ulam in 1941 \cite{kel1, kel2, ula1}, asserts that, when $n >2$, the converse also holds, i.e., $G$ is isomorphic to $H$ if and only if $G$ has the same collection of $n$ unlabelled cards as $H$. However, despite the efforts of many graph theorists, the status of the sufficiency of the condition remains unresolved. Surveys on the reconstruction problem can be found in \cite{bon2} \cite{b&h1} \cite{lau1}.

Since the conjecture remains unresolved, attention has focused on related reconstruction problems. One such area is proving that certain classes of graphs are {\em reconstructible} (i.e., that the conjecture is true when $G$ and $H$ belong to that class of graphs), or even {\em recognisable} (i.e., that membership of the class can be determined from the deck). Many classes of graphs have been shown to be reconstructible, including {\em trees} by Kelly \cite {kel2} and also Bondy \cite {bon1}, {\em disconnected graphs} by Greenwell and Hemminger \cite{g&h1}, and also Manvel \cite{man1}, {\em unicyclic graphs} by Manvel \cite{man2}, and {\em maximal planar graphs} by Lauri \cite {lau2}.

Another area of research has been to consider how many cards are required to reconstruct a graph - either its existential (ally) or universal (adversary) {\em reconstruction number} (see \cite{bbf1} \cite {myr1}) - or even just to recognise that it is a member of a particular class. An equivalent approach to finding universal reconstruction numbers is to consider the {\em maximum number of common cards} of two graphs. A {\em common card} of $G$ and $H$ is any card in the multiset intersection $\mathcal D(G) \cap \mathcal D(H)$, and the number of common cards of $G$ and $H$, denoted by $b(G,\, H)$, is the cardinality of this multiset intersection. The Reconstruction Conjecture can then be reformulated as follows: if $G$ and $H$ are not isomorphic then $b(G, \, H) < n$ when $n >2$.

Until a few years ago, there were no known families of pairs of non-isomorphic graphs that had $b(G,\, H) > \f{n}{2}+\f{1}{8}(3+\sqrt{8n+9})$. However, Bowler, Brown and Fenner \cite{bbf1} showed that there are, in fact, several infinite families of pairs of non-isomorphic graphs $G$ and $H$ with $b(G,\, H)=2\left\lfloor \f{n-1}{3} \right\rfloor$. Moreover, they conjectured that $b(G,\, H)$ is bounded above by $\f{2(n-1)}{3}$ for large enough $n$. Results for small graphs, i.e., for $n \leq 11$, have been provided by Baldwin \cite{bal1}, McMullen \cite{mcm1} and Rivshin \cite{riv1}.

In a subsequent paper \cite{bbf2}, Bowler, Brown, Fenner and Myrvold showed that if $G$ is disconnected and $H$ is connected then $b(G,\, H) \leq \left\lfloor \f{n}{2} \right\rfloor+1$, i.e., the connectedness of a graph can be recognised from any $\left\lfloor \f{n}{2} \right\rfloor+2$ of its cards. They also characterised all pairs of graphs that attain this bound  (most of these infinite families can also be found in \cite{bbf1}).

A similar recognition question is to ask how many cards are required to recognise whether a graph is a tree. Since $H$ is a tree in some of the families in \cite{bbf2}, it follows that at least $\left\lfloor \f{n}{2} \right\rfloor+2$ cards may be required. To fully answer the question, it is necessary to determine how many cards are required to distinguish a non-tree $G$ from a tree $H$ when it is known that $G$ is connected. It is easy to show that, in this case, $b(G,\,H) \le 2$ when $G$ contains more than one cycle. Furthermore, if $G-v \cong H-t$, for some $t \in V(H)$, then $v$ must be lie on every cycle of $G$. We can therefore restrict our attention to graphs where $G$ is unicyclic and the length of its cycle is reasonably large. For the maximum value of $b(G,\,H)$ to be attained, it has been conjectured that $G$ must be a {\em sunshine} graph (a connected graph for which $\mathrm{skel}(G)$ is a cycle) and $H$ must be a {\em caterpillar} (a connected graph for which $\mathrm{skel}(H)$ is a path); see \cite{bbf1}  \cite{fra1}. Support for restricting our attention to sunshine graphs and caterpillars is the observation that $b(G,\,H) \le 6$ when either $G$ is a sunshine graph and $H$ is a non-caterpillar tree, or when $G$ is a non-sunshine unicyclic graph and $H$ is a caterpillar.

In \cite{bro1}, Brown gave an intricate proof that, for large $n$, the number of common cards between a sunshine graph and a caterpillar of order $n$ is at most $2\left\lfloor\f{n+1}{5}\right\rfloor$ and, moreover, that this bound is only attained by a unique family of pairs of graphs for which $n \equiv 4 \pmod 5$. In this paper, we prove this result using {\em supercards}, a new approach to the study of the maximum number of common cards that we introduced in \cite{bf1}.

A {\em supercard} of $G$ and $H$ is a graph $G^+$, the deck of which contains at least one card isomorphic to $G$ and at least one card  isomorphic to $H$. In \cite{bf1}, we showed the existence of certain subsets of $S_{V(G^+)}$ of cardinality $b(G,\, H)$, the elements of which are in one-to-one correspondence with the elements of $\mathcal D(G) \cap \mathcal D(H)$. We called these subsets {\em maximum saturating sets}. We further showed that, for many families of pairs of graphs with a large number of common cards, there exist maximum saturating sets containing a large number of automorphisms of $G^+$.

We use supercards to investigate the case when $G$ is a sunshine graph and $H$ is a caterpillar. We show that, when $b(G,\, H) \ge 2\left\lfloor \sqrt{2n+1} \right\rfloor + 4$, we may construct supercards $G^+$ of $G$ and $H$, and suitable maximum saturating sets for which almost all of the elements are automorphisms of $G^+$. We further show that, in each case, the subset of elements that are automorphisms forms a group isomorphic to either a cyclic or dihedral group. We present several families of sunshine-caterpillar pairs that have such supercards and a large number of common cards. Finally, we show that $b(G,\, H)  \le \f{2(n+1)}{5}$ when $n \ge 62$, and we exhibit the unique family of pairs of graphs that attain this bound. In this case, the maximum saturating set is a subgroup of $\mathrm{Aut}(G^+)$ isomorphic to a dihedral group.
\section{\normalsize Supercards and Maximum Saturating Sets}\label{supercards section}\vspace{-2.5mm}
We now recall the main definitions and results in the theory of supercards. Detailed explanations and proofs of the results can be found in \cite{bf1}.
\begin{lemma}\label{iso cards and matching nhoods lemma}
{\rm Let $G$ and $H$ be graphs, and let $\gamma$ be a bijection from $V(G)$ to $V(H)$. Suppose that there is some vertex $v$ of $G$ such that $\gamma(G-v)=H-\gamma(v)$ and $\gamma(N_G(v))=N_H(\gamma(v))$. Then $\gamma(G)=H$.
\hfill\eot}
\end{lemma}
\begin{cor}\label{iso cards and matching nhoods aut cor}
{\rm Let $G$ be a graph and let $v \in V(G)$. Suppose that $\gamma \in S_{V(G)}$. Then $\gamma \in \mathrm{Aut}(G)$ if and only if 
\vspace{-2.5mm}
\[
\gamma(G-v)=G-\gamma(v) \qquad \text{and}  \qquad \gamma(N_G(v))=N_G(\gamma(v)).\vspace{-7mm}
\]
\hfill\eot}
\end{cor}
\begin{defn}\label{supercard definition 1}
{\rm A {\em supercard} of $G$ is any graph of order $n+1$ whose deck contains a card isomorphic to $G$.}
\end{defn}
\begin{defn}\label{supercard definition 2}
{\rm A {\em common supercard} of $G$ and $H$ is any graph that is a supercard of both $G$ and $H$, i.e., a graph whose deck contains some card $\widehat{G}$ isomorphic to $G$ and another card $\widehat{H}$ isomorphic to $H$. For brevity, we refer to such graphs as {\em supercards of $G$ and $H$}.}
\end{defn}
\begin{lemma}\label{supercard construction lemma 1}
{\rm There exists a graph $G^+$ that is a supercard of $G$ and $H$ if and only if $b(G,\, H) \geq 1$.

\vspace{-2.5mm}
\hfill\eot}
\end{lemma}
It is easy to verify that, if $v \in V(G)$, $t \in V(H)$ and $\gamma$ is an isomorphism such that $\gamma(G-v)=H-t$, then the graph $G^\ast$ defined by
\begin{eqnarray}
V(G^\ast)&=&V(G) \cup \left\{w\right\},\nonumber\\
E(G^\ast)&=&E(G) \cup \left\{xw \mid x \in V(G-v) \mbox{ and } \gamma(x)t \in E(H) \right\}, \label{supercard eqn 1}
\end{eqnarray}
for some $w \not\in V(G) \cup V(H)$, is a supercard of $G$ and $H$. In this case, $\widehat{G}=G^+-w = G$ and $\widehat{H}=G^+-v \cong H$ (see Lemma $3.3$ in \cite{bf1}.)

For the rest of this section we assume that $b(G,\, H) \ge 1$. We then let $G^+$ be a supercard of $G$ and $H$, and let $v$ and $w$ be vertices of $G^+$ such that $\widehat{G}=G^+-w \cong G$ and $\widehat{H}=G^+-v \cong H$. 
\begin{defn}\label{gamma-active pairs}
{\rm The set of {\em active permutations} of $G^+$ with respect to $v$ and $w$, denoted by $B_{vw}(G^+)$, is the subset of $S_{V(G^+)}$ defined by}
\begin{eqnarray}
B_{vw}(G^+)&=&\{\lambda \in S_{V(G^+)} \mid \lambda((G^+-w)-\lambda^{-1}(v))=(G^+-v)-\lambda(w) \}\\\label{B(G_hat_gamma)}
&=&\{\lambda \in S_{V(G^+)} \mid \lambda(\widehat{G}-\lambda^{-1}(v))=\widehat{H}-\lambda(w) \}. \nonumber
\end{eqnarray}
\end{defn}
We note that $1_{V(G^+)} \in B_{vw}(G^+)$, and that if $\lambda \in B_{vw}(G^+)$ then $\lambda(w) \ne v$, since $G$ and $H$ are not isomorphic.
\begin{defn}\label{B(G) defn}
{\rm $B^G_{vw}(G^+)$ is the subset of $B_{vw}(G^+)$ defined by}
\begin{equation}
B^G_{vw}(G^+)=\{\lambda \in B_{vw}(G^+) \mid \lambda(\widehat{G})=G^+-\lambda(w)\}.\label{B_G defn}
\end{equation}
\end{defn}
\begin{defn}\label{supercard number of common cards}
{\rm A {\em maximum saturating set} of $B_{vw}(G^+)$ is a subset $X \subseteq B_{vw}(G^+)$ that satisfies the following three properties:\vspace{-4mm}
\begin{enumerate}\lettnum
\item $1_{V(G^+)} \in X$;\vspace{-2.5mm}
\item if $\lambda$ and $\pi$ are distinct permutations in $X$ then $\lambda^{-1}(v) \ne \pi^{-1}(v)$ and $\lambda(w) \ne \pi(w)$;\vspace{-2.5mm}
\item there is no $\sigma$ in $B_{vw}(G^+)\setminus X$ such that $X \cup \left\{\sigma\right\}$ satisfies (b).\vspace{-2.5mm}
\end{enumerate}
}
\end{defn}
We note that, for any pair of distinct permutations $\lambda$ and $\pi$ in $X$, (b) guarantees that $G^+ - \lambda^{-1}(v) \ne G^+ - \pi^{-1}(v)$ and $G^+ - \lambda(w) \ne G^+ - \pi(w)$, although either pair of graphs could be isomorphic.

Although condition (c) only ensures that $X$ is {\em maximal} with respect to (a) and (b), it follows from Theorem \ref{max sat set lemma} below that all maximum saturating sets have the same cardinality. This implies that such sets are in fact of {\em maximum} cardinality with respect to (a) and (b), and justifies our terminology in Definition \ref{supercard number of common cards}.
\begin{defn}\label{B_G defn}
{\rm Let $X$ be a maximum saturating set of $B_{vw}(G^+)$. Then $X_G=X \cap B^G_{vw}(G^+)$.

\vspace{-2.5mm}
\hfill\eot}
\end{defn}
\begin{defn}\label{optimum max set defn}
{\rm A $G^+${\em-optimum saturating set} of $B_{vw}(G^+)$ is a maximum saturating set $X$ of $B_{vw}(G^+)$ for which $|X_G|$ takes its maximum possible value. We define $\chi(G^+)=|X_G|$ for any $G^+$-optimum saturating set $X$.
\hfill\eot}
\end{defn}
\begin{thm}\label{max sat set lemma}
{\rm Let $Y \subseteq B_{vw}(G^{+})$ satisfy properties (a) and (b) of Definition \ref{supercard number of common cards}. 
\begin{enumerate}\lettnum\vspace{-4mm}
\item If $Y$ is not a maximum saturating set of $B_{vw}(G^+)$ then $|Y| < b(G,\, H)$.\vspace{-2.5mm}
\item If $|Y| < b(G,\, H)$ then there is a maximum saturating set $X$ such that $Y \subset X$ (so $Y$ is not a maximum saturating set).\vspace{-2.5mm}
\item $Y$ is a maximum saturating set of $B_{vw}(G^+)$ if and only if $|Y|=b(G,\, H)$.\vspace{-2.5mm}
\end{enumerate}
\vspace{-8mm}
\hfill\eot}
\end{thm}
We make frequent use of the fact that every maximum saturating set of $B_{vw}(G^+)$ has cardinality $b(G,\, H)$ without explicitly quoting this theorem.
\section{\normalsize Sunshine graphs and caterpillars}\label{caterpillars and sunshine graphs}\vspace{-2.5mm}
We recall that $\mathrm {skel}(G)$, the skeleton of the graph $G$, is the graph $G-X$, where $X$ is the set of leaves and isolated vertices of $G$. A {\em sunshine graph} is a connected graph whose skeleton is a cycle and a {\em caterpillar} is a connected graph whose skeleton is a path. Clearly, all caterpillars are trees and all sunshine graphs are unicyclic. We shall use supercards to investigate the number of common cards between pairs of such graphs.

We denote the {\em diameter}, i.e., the length of a longest path, of a connected graph $G$ by $\mathrm {diam}(G)$. If $G$ consists of a tree $T$ plus a collection of isolated vertices, then we define $\mathrm {diam}(G)=\mathrm {diam}(T)$. A leaf at the end of any longest path in a graph is called a {\em peripheral} leaf. For any vertex $v$ of $G$, we denote the number of vertices in $N_G(v)$ of degree $2$ by $\tau_G(v)$. A $d${\em-leaf} of $G$  is a leaf $w$ (in a component of order at least three) for which $\tau_G(w)=0$, i.e., the degree of its neighbour is at least $3$. 

We are interested in pairs of sunshine graphs $U$ and caterpillars $T$ with a large number of common cards relative to their order $n$. We therefore assume, {\em for the rest of this paper}, that all the pairs $U$ and $T$ that we consider have $b(U,\, T) \ge 5$. By inspection, it is easy to show that, for such pairs, the unique cycle of $U$ is of length at least $6$, and that $U$ has at least one leaf. 

We use the following conventions for any sunshine graph $S$ with skeleton $x_0x_1 \ldots x_{c-1}x_0$: for any integer $k$, we interpret $x_k$ to be the vertex $x_i$, where $i$ is the unique integer such that $0 \leq i \leq c-1$ and $k \equiv i \pmod c$; in addition, for $0 \leq b < a \leq c-1$, we abbreviate the path $x_ax_{a+1} \ldots x_{c-1}x_0x_1 \ldots x_b$ on $\mathrm {skel}(S)$ to $x_ax_{a+1} \ldots x_b$.

We make frequent use of the following easy result concerning the cards of a sunshine graph.
\begin{lemma}\label{sunshine lemma 1}
{\rm Let $S$ be a sunshine graph with skeleton $x_0x_1 \ldots x_{c-1}x_0$ and let $x_i$ be in $V(\mathrm{skel}(S))$. Then $S-x_i$ consists of a caterpillar $Q$ of diameter $c-\tau_S(x_i)$, together with $d_S(x_i)-2$ isolated vertices. In addition:
\begin{enumerate}\lettnum\vspace{-4mm}
\item $\mathrm {skel}(Q)$ is $x_{i+1} x_{i+2} \ldots x_{i-1}$ if and only if $d_S(x_{i+1}) \geq 3$ and $d_S(x_{i-1}) \geq 3$, i.e., $\tau_S(x_i)=0$;\vspace{-2.5mm}
\item $\mathrm {skel}(Q)$ is $x_{i+1} x_{i+2} \ldots x_{i-2}$ if and only if $d_S(x_{i+1})\geq 3$ and $d_S(x_{i-1}) =2$, i.e., $\tau_S(x_i)=1$;\vspace{-2.5mm}
\item $\mathrm {skel}(Q)$ is $x_{i+2} x_{i+3} \ldots x_{i-1}$ if and only if $d_S(x_{i+1})=2$ and $d_S(x_{i-1}) \geq 3$, i.e., $\tau_S(x_i)=1$;\vspace{-2.5mm}
\item $\mathrm {skel}(Q)$ is $x_{i+2} x_{i+3} \ldots x_{i-2}$ if and only if $d_S(x_{i+1})=d_S(x_{i-1})=2$, i.e., $\tau_S(x_i)=2$.\vspace{-2.5mm}
\end{enumerate}
Moreover, $x_{i+1}$ is a peripheral leaf of $Q$ adjacent to $x_{i+2}$ when $d_S(x_{i+1}) =2$, i.e., in cases (c) and (d). Similarly, $x_{i-1}$ is a peripheral leaf of $Q$ adjacent to $x_{i-2}$ when $d_S(x_{i-1}) =2$, i.e., in cases (b) and (d).

\textit{Proof} Since $x_{i+1}$ and $x_{i-1}$ are the only non-leaves in $N_S(x_i)$ and $x_i$ is on the unique cycle of $S$, clearly $S-x_i$ consists of a tree $Q$ together with $d_S(x_i)-2$ isolated vertices. Moreover, since $\mathrm{skel}(S)$ is a cycle, $\mathrm{skel}(Q)$ is a path, so $Q$ is a caterpillar. It is easy to see that $x_{i+1}$ is a peripheral leaf of $Q$ when $d_S(x_{i+1})=2$; otherwise $x_{i+1}$ is one of the two vertices of $\mathrm{skel}(Q)$ that is adjacent to a peripheral leaf. A similar observation holds for $x_{i-1}$. Cases (a) to (d) then follow immediately. Finally, it follows from (a) to (d) that $\mathrm {diam}(Q)=c-\tau_S(x_i)$ as $\mathrm {diam}(Q)=|V(\mathrm {skel}(Q))|+1$.
\hfill\eot}
\end{lemma}
We also make the following easy observation about the possible cards of a caterpillar.
\begin{lemma}\label{cat lemma 1}
{\rm Let $T$ be a caterpillar with skeleton $y_1y_2 \ldots y_p$.
\begin{enumerate}\lettnum\vspace{-4mm}
\item $T-y_i$ contains precisely one non-trivial component if and only if $i \in \{1,\, p\}$.\vspace{-2.5mm}
\item If $t$ is a leaf of $T$ that is not a $d$-leaf then $t$ is a peripheral leaf, and is adjacent to either $y_1$ or $y_p$. Moreover, $t$ is the only leaf that is adjacent to $y_1$ or $y_p$, respectively.\vspace{-2.5mm}
\item There exist at most two leaves of $T$ that are not $d$-leaves.\vspace{-2.5mm}
\end{enumerate}
\textit{Proof} The results follow immediately by considering the structure of $T$.
\hfill\eot}
\end{lemma}
These two results yield the following important lemma.
\begin{lemma}\label{sun cat lemma new 2}
{\rm Let $U$ be a sunshine graph and $T$ be a caterpillar. Then there exists a supercard of $U$ and $T$ that is a sunshine graph.

\textit{Proof} Since $U$ is a sunshine graph, clearly no card of $U$, and hence no common card of $U$ and $T$, can contain more than one non-trivial component. So, since $b(U,\, T) \ge 5$, it follows from Lemma \ref{cat lemma 1} that there exists a vertex $v$ of $U$, a $d$-leaf $t$ of $T$, and an isomorphism $\gamma$ such that $\gamma(U-v)=T-t$.

Let $s$ be the unique vertex of $T$ adjacent to $t$, and let $U^\ast$ be the supercard of $U$ and $T$ constructed as in equation (\ref{supercard eqn 1}). Now, since $t$ is a leaf of $T$, clearly $\gamma^{-1}(s)$ is the only vertex of $U^\ast$ adjacent to $w$. Moreover, since $t$ is a $d$-leaf of $T$, it is easy to see that $d_U(\gamma^{-1}(s)) \ge 2$, so $\gamma^{-1}(s)$ is on the unique cycle of $U$. Thus $U^\ast$ must be a sunshine graph. 
\hfill\eot}
\end{lemma}
Let $U$ be a sunshine graph and $T$ a caterpillar. By Lemma \ref{sun cat lemma new 2}, there exists some supercard $U^+$ of $U$ and $T$ that is a sunshine graph. So, for any such supercard $U^+$, let $w$ and $v$ be vertices of $U^+$ such that $\widehat{U}=U^+-w \cong U$ and $\widehat{T}=U^+-v \cong T$. Clearly, $w$ is a leaf of $U^+$ as $U$ is a sunshine graph. In addition, since $T$ is a tree, $v$ must be on the cycle of $U^+$ and $d_{U^+}(v)=2$. We therefore label the skeleton of $U^+$ as $x_0x_1 \ldots x_{c-1}x_0$ (so the cycle is of length $c \ge 6$), where $x_0$ is adjacent to $w$, $x_\nu$ is $v$, for some $\nu$, $1 \le \nu \le c-1$, and $d_{U^+}(x_{\nu-1}) \ge d_{U^+}(x_{\nu+1})$. We further arbitrarily label all the leaves of $U^+$, so that each distinct leaf adjacent to $x_i$ is labelled $x_i^j$ for some unique $j$, $1 \le j \le d_{U^+}(x_i)-2$, where $w$ is labelled $x_0^1$. Our supercard $U^+$ thus satisfies
\begin{equation}
V(U^+)=V(U) \cup \{w\} \quad \text{and }\quad E(U^+)=E(U) \cup \{x_0w\}, \label{sunshine supercard defn}
\end{equation}
assuming the above restrictions on $\mathrm {skel}(U^+)$.

{\em For the rest of this section}, we assume that $U^+$ is the supercard of $U$ and $T$ specified in (\ref{sunshine supercard defn}). Clearly, any supercard of $U$ and $T$ is also a supercard of any pair of graphs isomorphic to $U$ and $T$, respectively. So, for {\em ease of notation}, we shall write $U=U^+-w$ instead $\widehat{U}=U^+-w$, and $T=U^+-x_\nu$, instead of $\widehat{T}=U^+-x_\nu$. We note that $\mathrm {skel}(U)=\mathrm {skel}(U^+)$. Applying Lemma \ref{sunshine lemma 1} to $U^+$ and $x_\nu$ yields the following result.
\begin{lemma}\label{supercard sunshine lemma new 1}
{\rm We have the following possibilities for the skeleton of $T$.
\begin{enumerate}\lettnum\vspace{-2.5mm}
\item If $\tau_{U^+}(x_\nu)=0$ then $\mathrm {skel}(T)$ is $x_{\nu+1} x_{\nu+2} \ldots x_{\nu-1}$.\vspace{-2.5mm}
\item If $\tau_{U^+}(x_\nu)=1$ then $\mathrm {skel}(T)$ is $x_{\nu+2} x_{\nu+3} \ldots x_{\nu-1}$.\vspace{-2.5mm}
\item If $\tau_{U^+}(x_\nu)=2$ then $\mathrm {skel}(T)$ is $x_{\nu+2} x_{\nu+3} \ldots x_{\nu-2}$.\vspace{-1mm}
\end{enumerate}
It follows that $x_{\nu+1}$ is a peripheral leaf of $T$ adjacent to $x_{\nu+2}$ when $\tau_{U^+}(x_\nu) \ge 1$, and $x_{\nu-1}$ is a peripheral leaf of $T$ adjacent to $x_{\nu-2}$ when $\tau_{U^+}(x_\nu) = 2$.

\textit{Proof} This follows immediately by Lemma \ref{sunshine lemma 1} with $S=U^+$ and $i=\nu$, noting that case (b) of that lemma cannot occur as $d_{U^+}(x_{\nu-1}) \ge d_{U^+}(x_{\nu+1})$.
\hfill\eot}
\end{lemma}
We recall that $B_{vw}(U^+)$ is the set of active permutations of $U^+$ with respect to $v$ and $w$, i.e., 
\begin{equation}
B_{vw}(U^+)=B_{x_\nu x_0^1}(U^+)=\{\lambda \in S_{V(U^+)} \mid \lambda(U-\lambda^{-1}(x_\nu))=T-\lambda(w) \}.\label{B(U_hat_gamma)}
\end{equation}
\begin{lemma}\label{supercard sunshine lemma 1a}
{\rm Let $\lambda \in B_{vw}(U^+)$. Then exactly one of the following holds:
\begin{enumerate}\lettnum\vspace{-4mm}
\item $\lambda(w)$ is a $d$-leaf of $T$, in which case $\mathrm {skel}(T-\lambda(w))=\mathrm {skel}(T)$ and $\mathrm{diam}(T-\lambda(w))=\mathrm{diam}(T)$;\vspace{-2.5mm}
\item $\lambda(w)$ is a peripheral leaf of $T$ that is not a $d$-leaf, in which case $\mathrm{diam}(T-\lambda(w)))=\mathrm{diam}(T)-1$;\vspace{-2.5mm}
\item $\lambda(w)$ is adjacent to a peripheral leaf of $T$, in which case $\mathrm{diam}(T-\lambda(w)) \le \mathrm{diam}(T)-1$.\vspace{-2.5mm}
\end{enumerate}
\textit{Proof} By Lemmas \ref{sunshine lemma 1} and \ref{cat lemma 1}(a), $\lambda(w)$ is either a leaf or adjacent to a peripheral leaf. Cases (a), (b) and (c) then follow easily by considering the structure of $T$.
\hfill\eot}
\end{lemma}
\begin{lemma}\label{supercard sunshine lemma 1}
{\rm Let $\lambda \in B_{vw}(U^+)$. Then $\lambda^{-1}(x_\nu) \in V(\mathrm{skel}(U))$, i.e., $\lambda^{-1}(x_\nu)$ is $x_\mu$ for some $\mu$, $0 \le \mu \le c-1$, and $d_U(x_\mu)=d_T(\lambda(w))+1$. In addition, $U-x_\mu$, and therefore also $T-\lambda(w)$, consists of a caterpillar of diameter $c-\tau_U(x_\mu)$, together with $d_U(x_\mu)-2$, equivalently $d_T(\lambda(w))-1$, isolated vertices.

\textit{Proof} Since $T$ is a tree, clearly $\lambda^{-1}(x_\nu) \in V(\mathrm{skel}(U))$, so $\lambda^{-1}(x_\nu)$ is $x_\mu$, for some $\mu$. In addition, $d_U(x_\mu)=d_T(\lambda(w))+1$ as $|E(U)|=|E(T)|+1$. Finally, by Lemma \ref{sunshine lemma 1}, $U-x_\mu$ consists of a caterpillar of diameter $c-\tau_U(x_\mu)$, together with $d_U(x_\mu)-2$ isolated vertices.
\hfill\eot}
\end{lemma}
From Lemma \ref{supercard sunshine lemma 1} and Theorem \ref{max sat set lemma}, it immediately follows that $b(U,\, T) \le c$.

For brevity, we frequently use Lemma \ref{supercard sunshine lemma 1} in the rest of this paper without explicit reference. Moreover, given any $\lambda \in B_{uv}(U^+)$, unless otherwise stated, we let $x_\mu=\lambda^{-1}(x_\nu)$. So, by Definition \ref{gamma-active pairs},
\begin{equation}
\lambda(U-x_\mu)=\lambda((U^+-w)-x_\mu)=(U^+-x_\nu)-\lambda(w)=T-\lambda(w).\label{sun cat B(G+) eqn}
\end{equation}
We write the skeletons of the caterpillars in $U-x_\mu$ and ${T}-\lambda(w)$ as
\begin{equation}
\mathrm{skel}(U-x_\mu):\, x_ax_{a+1} \ldots x_b \qquad \text{and} \qquad \mathrm{skel}(T-\lambda(w)):\, x_rx_{r+1} \ldots x_s,\label{cat skels}
\end{equation}
respectively. On applying Lemma \ref{sunshine lemma 1} with $S=U$ and $i=\mu$, we have the following result.
\begin{cor}\label{supercard sunshine lemma 2}
{\rm Let $\mathrm{skel}(U-x_\mu)$ be as in (\ref{cat skels}). Then $x_a \in \{x_{\mu+1},\, x_{\mu+2}\}$ and $x_b \in \{x_{\mu-1},\, x_{\mu-2}\}$.

\vspace{-1mm}
\hfill\eot}
\end{cor}
\vspace{-5mm}
Since $c \ge 6$, it immediately follows that $a \ne b$.
\begin{lemma}\label{supercard sunshine lemma 4}
{\rm Let $\lambda \in B_{vw}(U^+)$. Then $b-a \equiv s-r \pmod c$. Moreover, either
\begin{enumerate}\lettnum\vspace{-2.5mm}
\item $\lambda(x_i)=x_{(r-a)+i}$ for all $x_i \in V(\mathrm{skel}(U-x_\mu))$, or\vspace{-2mm}
\item $\lambda(x_i)=x_{(s+a)-i}$ for all $x_i \in V(\mathrm{skel}(U-x_\mu))$.\vspace{-2mm}
\end{enumerate}
We note that $\lambda(x_b)=x_s$ in (a) and $\lambda(x_b)=x_r$ in (b).

\textit{Proof} $\lambda$ maps the skeleton of $U-x_\mu$ onto the skeleton of $T-\lambda(w)$. So $b-a \equiv s-r \pmod c$, and either $\lambda(x_a)=x_r$ and $\lambda(x_b)=x_s$, or $\lambda(x_a)=x_s$ and $\lambda(x_b)=x_r$. It is then easy to see that either (a) or (b) must hold.
\hfill\eot}
\end{lemma} 
\begin{lemma}\lettnum\label{supercard sunshine lemma 5}
{\rm Let $\lambda \in B_{vw}(U^+)$. Then
\begin{equation}
c-\tau_{U^+}(x_\nu) = \mathrm{diam}(T) \ge \mathrm{diam}(T-\lambda(w))=\mathrm{diam}(U-x_\mu)=c-\tau_U(x_\mu) \ge c-2.\label{diam eqn}
\end{equation}
\textit{Proof} This follows easily by applying Lemma \ref{sunshine lemma 1} to $U^+$ and $x_\nu$, and then $U$ and $x_\mu$.
\hfill\eot}
\end{lemma} 
\begin{lemma}\lettnum\label{supercard sunshine lemma 5a_2}
{\rm Let $\lambda \in B_{vw}(U^+)$ be such that $\lambda(w)$ is a $d$-leaf of $T$. Suppose that $\lambda(x_{\mu+2}) \not\in \left\{x_{\nu-2},\, x_{\nu+2} \right\}$. Then $\tau_U(x_\mu)=\tau_{U^+}(x_\nu)=1$. Moreover, $\mathrm{skel}(U-x_\mu)$ is either
\begin{enumerate}\lettnum\vspace{-2.5mm}
\item $x_{\mu+1}x_{\mu+2} \ldots x_{\mu-2}$, in which case $\lambda(x_i)= x_{(\nu-\mu+1)+i}$ for all $x_i$ in $V(\mathrm{skel}(U-x_\mu))$, or\vspace{-2.5mm}
\item $x_{\mu+2}x_{\mu+3} \ldots x_{\mu-1}$, in which case $\lambda(x_i)= x_{(\nu+\mu+1)-i}$ for all $x_i$ in $V(\mathrm{skel}(U-x_\mu))$.\vspace{-2.5mm}
\end{enumerate}
\textit{Proof} $\mathrm{skel}(T-\lambda(w))=\mathrm{skel}(T)$ by Lemma \ref{supercard sunshine lemma 1a}(a). So it follows from Corollary \ref{supercard sunshine lemma new 1} that the possible skeletons of $T-\lambda(w)$ are determined by $\tau_{U^+}(x_\nu)$. In addition, $\tau_U(x_\mu)=\tau_{U^+}(x_\nu)$ by (\ref{diam eqn}). It then follows from Lemma \ref{sunshine lemma 1} that the possible skeletons of $U-x_\mu$ are also determined by $\tau_{U^+}(x_\nu)$. Using Lemma \ref{supercard sunshine lemma 4}, it is now straightforward to determine all the possibilities for $\lambda(x_{\mu+2})$, for each of the three values of $\tau_{U^+}(x_\nu)$.

It is easy to see that $\lambda(x_{\mu+2}) \in \{x_{\nu-2},\, x_{\nu+2} \}$ when $\tau_{U^+}(x_\nu) \ne 1$. Since this excluded by assumption, it follows that $\tau_U(x_\nu)=\tau_{U^+}(x_\mu)=1$. In this case, $x_r$ is $x_{\nu+2}$ and $x_s$ is $x_{\nu-1}$ by Corollary \ref{supercard sunshine lemma new 1}(b). So, since $\lambda(x_{\mu+2}) \not\in \left\{x_{\nu-2},\, x_{\nu+2} \right\}$, it is straightforward to show that either (i) $x_a$ is $x_{\mu+1}$, $x_b$ is $x_{\mu-2}$ and Lemma \ref{supercard sunshine lemma 4}(a) holds, or (ii) $x_a$ is $x_{\mu+2}$, $x_b$ is $x_{\mu-1}$ and Lemma \ref{supercard sunshine lemma 4}(b) holds. In case (i), we have $r-a \equiv \nu-\mu+1 \pmod c$, and in case (ii),  we have $s+a \equiv \nu+\mu+1 \pmod c$. Cases (a) and (b) of the lemma then immediately follow from Lemma \ref{supercard sunshine lemma 4}.
\hfill\eot}
\end{lemma}
\begin{defn}\lettnum\label{supercard sunshine defn X_tilde}
{\rm We define $\widetilde{B}_{vw}(U^+)$ to be the subset of $ B_{vw}(U^+)$ containing those permuations $\lambda$ such that $\lambda(w)$ is a leaf of $U^+$ and a $d$-leaf of $T$. We further define $\widetilde{X} = X \cap \widetilde{B}_{vw}(U^+)$ for any maximum saturating set $X$ of $B_{vw}(U^+)$.
\hfill\eot}
\end{defn}
\begin{lemma}\label{supercard sunshine lemma 4b}
{\rm Let $\lambda \in B_{vw}(U^+) \setminus \widetilde{B}_{vw}(U^+)$.
\begin{enumerate}\lettnum\vspace{-4mm}
\item If $\tau_{U^+}(x_\nu)=0$ then $\lambda(w)$ is not a $d$-leaf of $T$. Moreover, 
\begin{enumerate}\vspace{-2.5mm}
\item[(i)] if $\lambda(w)$ is a leaf of $T$ then either $\lambda(w)=x_{\nu+1}^1$ and $d_{U^+}(x_{\nu+1})=3$, or $\lambda(w)=x_{\nu-1}^1$ and $d_{U^+}(x_{\nu-1})=3$;\vspace{-1mm}
\item[(ii)] if $\lambda(w)$ is a cut-vertex of $T$ then $\lambda(w) \in \{x_{\nu+1}, x_{\nu-1}\}$.\vspace{-2mm}
\end{enumerate}
\item If $\tau_{U^+}(x_\nu)=1$ then 
\begin{enumerate}\vspace{-2.5mm}
\item[(i)] if $\lambda(w)$ is a $d$-leaf of $T$ then $\lambda(w)=x_{\nu+1}$ and $d_{U^+}(x_{\nu+2}) \ge 3$;\vspace{-1mm}
\item[(ii)] if $\lambda(w)$ is a leaf but not a $d$-leaf of $T$ then either $\lambda(w)=x_{\nu+1}$ and $d_{U^+}(x_{\nu+2})=2$, or $\lambda(w)=x_{\nu-1}^1$ and $d_{U^+}(x_{\nu-1})=3$;\vspace{-1mm}
\item[(iii)] if $\lambda(w)$ is a cut-vertex of $T$ then $\lambda(w) \in \{x_{\nu+2}, x_{\nu-1}\}$ \vspace{-2mm}
\end{enumerate}
\item If $\tau_{U^+}(x_\nu)=2$ then $\lambda(w)$ is a $d$-leaf of $T$ and $\lambda(w) \in \{x_{\nu+1},\, x_{\nu-1}\}$.\vspace{-2mm}
\end{enumerate}
\textit{Proof} Since $\lambda \not\in \widetilde{B}_{vw}(U^+)$, either $\lambda(w)$ is a $d$-leaf of $T$ that is not a leaf of $U^+$, or it is not a $d$-leaf of $T$, in which case (b) or (c) of Lemma \ref{supercard sunshine lemma 1a} must hold.

(a) Suppose that $\tau_{U^+}(x_\nu)=0$. Then $\mathrm{skel}(T)$ is given by Corollary \ref{supercard sunshine lemma new 1}(a), so every leaf of $T$ is a leaf of $U^+$. It immediately follows that $\lambda(w)$ cannot be a $d$-leaf of $T$. Now, if case (b) of Lemma \ref{supercard sunshine lemma 1a} holds then it is easy to see that either $\lambda(w)=x^1_{\nu+1}$ and $d_{U^+}(x_{\nu+1})=3$, or $\lambda(w)=x_{\nu-1}^1$ and $d_{U^+}(x_{\nu-1})=3$. On the other hand, if case (c) holds then $\lambda(w)$ is either $x_{\nu+1}$ or $x_{\nu-1}$.

(b) Suppose that $\tau_{U^+}(x_\nu)=1$. Then $\mathrm{skel}(T)$ is given by Corollary \ref{supercard sunshine lemma new 1}(b), so the only possible $d$-leaf of $T$ that is not a leaf of $U^+$ is $x_{\nu+1}$; in this case, clearly $d_{U^+}(x_{\nu+2}) \ge 3$. Now, if case (b) of Lemma  \ref{supercard sunshine lemma 1a} holds then it is easy to see that either  $\lambda(w)=x_{\nu+1}$ and $d_{U^+}(x_{\nu+2})=2$, or $\lambda(w)=x_{\nu-1}^1$ and $d_{U^+}(x_{\nu-1})=3$. On the other hand, if case (c) holds then $\lambda(w)$ is either $x_{\nu+2}$ or $x_{\nu-1}$.

(c) Suppose that $\tau_{U^+}(x_\nu)=2$. Then $\mathrm{skel}(T)$ is given by Corollary \ref{supercard sunshine lemma new 1}(c), so the only possible $d$-leaves of $T$ that are not leaves of $U^+$ are $x_{\nu+1}$ and $x_{\nu-1}$. Since equality holds throughout (\ref{diam eqn}), clearly $\mathrm{diam}(T-\lambda(w)) = \mathrm{diam}(T)$ and, therefore, neither case (b) nor case (c) of Lemma \ref{supercard sunshine lemma 1a} holds.
\hfill\eot}
\end{lemma}
\begin{cor}\label{supercard sunshine cor 4c}
{\rm Let $X$ be a maximum saturating set of $B_{vw}(U^+)$.
\begin{enumerate}\lettnum\vspace{-4mm}
\item If $\tau_{U^+}(x_\nu)=0$ then $|X \setminus \widetilde{X}| \le 4$.\vspace{-2mm}
\item If $\tau_{U^+}(x_\nu)=1$ then $|X \setminus \widetilde{X}| \le 4$.\vspace{-2mm}
\item If $\tau_{U^+}(x_\nu)=2$ then $|X \setminus \widetilde{X}| \le 2$.\vspace{-2mm}
\end{enumerate}
\textit{Proof} This follows immediately from Lemma \ref{supercard sunshine lemma 4b} and part (b) of Definition \ref{supercard number of common cards}.
\hfill\eot}
\end{cor}
We recall from Definition \ref{optimum max set defn} that a $U^+$-optimum saturating set $X$ of $B_{vw}(U^+)$ is a maximum saturating set of $B_{vw}(U^+)$ such that $|X_U|= \chi(U^+)$.
\begin{cor}\label{supercard sunshine cor 4d}
{\rm If there exists a $U^+$-optimum saturating set $X$ of $B_{vw}(U^+)$ such that $\widetilde{X} \subseteq X_U$ then $b(U,\, T) \le  \chi(U^+) + 4$.
\hfill\eot}
\end{cor}
We now consider the permuations in $B^U_{vw}(U^+)$, i.e., those permutations $\lambda$ in $B_{vw}(U^+)$ such that $\lambda(U)=U^+-\lambda(w)$.
\begin{lemma}\label{supercard sunshine lemma 6}
{\rm Let $\lambda \in B^U_{vw}(U^+)$. Then $\lambda(w)$ is a leaf of $U^+$ adjacent to a vertex of degree $d_{U^+}(x_0)$, $d_U(x_\mu)=2$, and either
\begin{enumerate}\lettnum\vspace{-4mm}
\item $\lambda(x_i)=x_{(\nu-\mu)+i}$ for all $x_i$ (a rotation), or\vspace{-2mm}
\item $\lambda(x_i)=x_{(\nu+\mu)-i}$ for all $x_i$ (a reflection).\vspace{-2mm}
\end{enumerate}
\textit{Proof} Since $\lambda(U)=U^+-\lambda(w)$, clearly $\lambda(w)$ must be a leaf of $U^+$ adjacent to a vertex of degree $d_{U^+}(x_0)$. Thus $d_U(x_\mu)=2$ by Lemma \ref{supercard sunshine lemma 1}. Now $\lambda$ must preserve the cycle structure of $U$. So, since $\lambda(x_\mu)=x_\nu$, it follows that either $\lambda(x_{\mu+i})=x_{\nu+i}$ and $\lambda(x_{\mu-i})=x_{\nu-i}$ for all $i$, or $\lambda(x_{\mu+i})=x_{\nu-i}$ and $\lambda(x_{\mu-i})=x_{\nu+i}$ for all $i$. Cases (a) and (b) then follow immediately.
\hfill\eot}
\end{lemma} 
The following lemma gives a methodology for replacing permutations in $\widetilde{B}_{vw}(U^+)$ by ``equivalent'' permutations in $B^U_{vw}(U^+)$.
\begin{lemma}\label{supercard sunshine lemma 7}
{\rm Let $\lambda \in \widetilde{B}_{vw}(U^+)$. Suppose that $\lambda(x_{\mu+2}) \in \left\{x_{\nu-2},\, x_{\nu+2} \right\}$. Then there exists $\widehat{\lambda} \in  B^U_{vw}(U^+)$ such that $\widehat{\lambda}^{-1}(x_\nu)=\lambda^{-1}(x_\nu)=x_\mu$ and $\widehat{\lambda}(w)=\lambda(w)$.

\textit{Proof} By Corollary \ref{supercard sunshine lemma 2}, $x_a \in \{x_{\mu+1},\, x_{\mu+2}\}$ and $x_b \in \{x_{\mu-1},\, x_{\mu-2}\}$; so $x_{\mu+2} \in  \{x_a,\, x_{a+1}\}$. Since $\lambda(w)$ is a $d$-leaf of $T$, it follows from Lemma \ref{supercard sunshine lemma 1a}(a) that $\mathrm {skel}(T-\lambda(w))=\mathrm {skel}(T)$. So $x_r \in \{x_{\nu+1},\, x_{\nu+2}\}$ and $x_s \in \{x_{\nu-1},\, x_{\nu-2}\}$ by Corollary \ref{supercard sunshine lemma new 1}. We now assume that $\lambda(x_{\mu+2}) =x_{\nu+2}$ and prove the result in this case. The case when $\lambda(x_{\mu+2}) =x_{\nu-2}$ can be proved in a similar manner.

Suppose that Lemma \ref{supercard sunshine lemma 4}(b) holds. Then $\{\lambda(x_a),\, \lambda(x_{a+1})\} \subseteq \{x_s,\, x_{s-1}\}  \subseteq \{x_{\nu-1}, \, x_{\nu-2}, \, x_{\nu-3}\}$. Since $\lambda(x_{\mu+2}) =x_{\nu+2}$ and $x_{\mu+2} \in  \{x_a,\, x_{a+1}\}$, this implies that $x_{\nu+2} \in \{x_{\nu-1}, \, x_{\nu-2}, \, x_{\nu-3}\}$. This is impossible as $c \ge 6$. Therefore Lemma \ref{supercard sunshine lemma 4}(a) must hold. Hence $r-a \equiv s-b \equiv \nu-\mu \pmod c$, and $\lambda(x_i)=x_{(\nu-\mu)+i}$ for all $x_i$ in $V({\rm skel}(U-x_\mu))$.

Let $\theta$ be the transposition of $V(T-\lambda(w))$ that swaps $\lambda(x_{\mu+1})$ and $x_{\nu+1}$. We show that $\theta \in \mathrm{Aut}(T-\lambda(w))$. If $\lambda(x_{\mu+1}) = x_{\nu+1}$ then $\theta$ is $1_{V(T-\lambda(w))}$, so there is nothing to prove. Suppose therefore that $\lambda(x_{\mu+1}) \ne x_{\nu+1}$. It is then easy to show that $x_{\mu+1} \not\in V(\mathrm{skel}(U-x_\mu))$, so $x_a$ is $x_{\mu+2}$ and $x_{\mu+1}$ is a leaf of $U-x_\mu$ adjacent to $x_a$. Hence, correspondingly, $\lambda(x_{\mu+1})$ is a leaf of $T-\lambda(w)$ adjacent to $x_r$. Furthermore, $x_r$ is $x_{\nu+2}$, so $x_{\nu+1}$ is also a leaf of $T-\lambda(w)$ adjacent to $x_r$. Thus $\theta$ swaps a pair of leaves adjacent to $x_{\nu+2}$, so $\theta \in \mathrm{Aut}(T-\lambda(w))$ in this case also.

By considering the vertices $x_b$ and $x_s$, it is easy to show that the transposition $\theta^\prime$ of $V(T-\lambda(w))$ that swaps $\lambda(x_{\mu-1})$ and $x_{\nu-1}$ is also in $\mathrm{Aut}(T-\lambda(w))$. Let us define $\widehat{\lambda} \in S_{V(U^+)}$ by $\widehat{\lambda}(x_\mu)=x_\nu$, $\widehat{\lambda}(w)=\lambda(w)$ and $\widehat{\lambda}(u)=\theta^\prime\theta\lambda(u)$ for all other vertices $u$ of $U^+$. Then
\[
\widehat{\lambda}(U-x_\mu)=\theta^\prime\theta\lambda(U-x_\mu)=\theta^\prime\theta(T-\lambda(w))=T-\lambda(w),
\]
so $\widehat{\lambda} \in B_{vw}(U^+)$. Moreover, since $\widehat{\lambda}(x_{\mu+1})=x_{\nu+1}$ and $\widehat{\lambda}(x_{\mu-1})=x_{\nu-1}$, it follows that $\widehat{\lambda}(N_U(x_\mu))=\{x_{\nu-1},\, x_{\nu+1} \}$ as $d_U(x_\mu)=2$. Hence $\widehat{\lambda}(U)=U^+-\lambda(w)$ by Lemma \ref{iso cards and matching nhoods lemma}, i.e., $\widehat{\lambda} \in  B^U_{vw}(U^+)$.
\hfill\eot}
\end{lemma}
\begin{lemma}\label{supercard sunshine lemma 7 new}
{\rm Let $X$ be a maximum saturating set of $B_{vw}(U^+)$ and let $L$ be the subset of $X$ defined by $L=\big \{\lambda \in \widetilde{X} \mid \lambda(x_{\mu+2}) \in \left\{x_{\nu-2},\, x_{\nu+2} \right\} \big \}$. Then
\begin{enumerate}\lettnum\vspace{-2.5mm}
\item $|L| \le \chi(U^+)$;\vspace{-2mm}
\item if $X$ is a $U^+$-optimum saturating set of $B_{vw}(U^+)$ then $X_U \cap \widetilde{X} = L$.\vspace{-2mm}
\end{enumerate}
\textit{Proof} For each $\lambda \in L$, we let $\widehat{\lambda}$ be the permutation in $B^U_{vw}(U^+)$ as defined in Lemma \ref{supercard sunshine lemma 7}. We then define $\widehat{L}$ to be the set of all such $\widehat{\lambda}$. Since $X$ is a maximum saturating set of $B_{vw}(U^+)$, clearly $X \setminus L$ and $\widehat{L}$ are disjoint, and the set $\widehat{X}$ defined by $\widehat{X}=(X \setminus L)  \cup \widehat{L}$ is also a maximum saturating set. Moreover $\widehat{X}_U = (X_U \setminus L) \cup \widehat{L}$. Therefore $|L| =|\widehat{L}| \le |\widehat{X}_U| \le \chi(U^+)$.

Now suppose that $X$ is $U^+$-optimum. If $\lambda \in X_U \cap \widetilde{X}$ then $\lambda \in L$ by Lemma \ref{supercard sunshine lemma 6}. Thus $X_U \cap \widetilde{X} \subseteq L$. So suppose that $L \not\subseteq X_U$. Then $|\widehat{X}_U|=|X_U \setminus L| + |\widehat{L}| > |X_U|$. This is impossible since $X$ is $U^+$-optimum. So $L \subseteq X_U$, and therefore $X_U \cap \widetilde{X} = L$.
\hfill\eot}
\end{lemma}
\begin{cor}\label{supercard sunshine cor 7a}
{\rm Let $X$ be a $U^+$-optimum saturating set of $B_{vw}(U^+)$. If there exists $\lambda \in  \widetilde{X} \setminus X_U$ then $\tau_{U^+}(x_\nu) = 1$ and $\{\lambda(x_{\mu+2}),\, \lambda(x_{\mu-2})\}  = \{x_{\nu+3},\, x_{\nu-1}\}$.

\textit{Proof} This follows easily from Lemma \ref{supercard sunshine lemma 7 new}(b) and Lemma \ref{supercard sunshine lemma 5a_2}.
\hfill\eot}
\end{cor}
It follows from this result that $\widetilde{X} \subseteq X_U$ when $\tau_{U^+}(x_\nu) \ne 1$.

If there exists a $U^+$-optimum saturating set $X$ of $B_{vw}(U^+)$ such that $\widetilde{X} \subseteq X_U$ then $b(U,\, T) \le \chi(U^+)+4$ by Corollary \ref{supercard sunshine cor 4d}. When there is no such set, we must construct another supercard of $U$ and $T$.

Let $\sigma \in \widetilde{B}_{vw}(U^+)$, and let $x_\xi=\sigma^{-1}(x_\nu)$, so $\sigma(U-x_\xi)=T- \sigma(w)$. Let $u$ be the unique vertex of $U$ adjacent to $\sigma(w)$. Since $\sigma(w)$ is a $d$-leaf of $T$, clearly $d_{T-\sigma(w)}(u) \ge 2$, and thus $d_{U-x_\xi}(\sigma^{-1}(u)) \ge 2$. Hence $\sigma^{-1}(u)$ is $x_\eta$ for some $\eta$, $0 \le \eta \le c-1$. We note that $x_\eta$ cannot be $x_\xi$.

We now define a new sunshine graph $U^+_\sigma$ constructed from $U^+$. First we delete from $U^+$ the edge $x_0w$ and add an additional edge $x_\eta w$. We then relabel the skeleton of $U^+$ as $z_0z_1 \ldots z_{c-1}z_0$, where $x_\eta$ is relabelled as $z_0$, $x_\xi$ as $z_\zeta$, and the other vertices on the cycle in the natural way, reversing the labelling around the cycle if necessary in order to ensure that $d_{U^+_\sigma}(z_{\zeta-1}) \ge d_{U^+_\sigma}(z_{\zeta+1})$. Finally, as before, we relabel all the leaves of $U^+$ so that each distinct leaf adjacent to $z_i$ is labelled $z_i^j$ for some unique $j$, $1 \le j \le d_{U^+_\sigma}(z_i)-2$, where $w$ is labelled $z_0^1$. We note this labelling is analogous to the original labelling of $U^+$.

Since $|V(U^+)|=|V(U^+_\sigma)|$, we may define a bijection $\theta$ from $V(U^+)$ to $V(U^+_\sigma)$ that encapsulates the relabelling described above, i.e., $\theta(x_i)=z_{i-\eta}$ if the order of the labels around the cycle did not need reversing, and $\theta(x_i)=z_{\eta-i}$ if it did. We also specify that $\theta(w)=w$, and that $\theta$ maps the remaining leaves of $U^+$ so that those adjacent to $x_i$ map to leaves adjacent to $\theta(x_i)$ for each $x_i$. We note that $\theta(x_\eta)=z_0$ and $\theta(x_\xi)=z_\zeta$.

Clearly, $d_{U^+_\sigma}(z_\zeta)=2$, as $d_U(x_\xi)=2$ and $x_\eta$ is not $x_\xi$. Let $U_\sigma=U^+_\sigma-w$ and $T_\sigma=U^+_\sigma-z_\zeta$. Then $U_\sigma^+$ is a supercard of $U_\sigma$ and $T_\sigma$ that satisfies
\begin{equation}
V(U^+_\sigma)=V(U_\sigma) \cup \{w\} \quad \text{and }\quad E(U^+_\sigma)=E(U_\sigma) \cup \{z_0w\}. \label{sunshine supercard 2 defn}
\end{equation}
For ease of notation, we write $B_{vw}(U^+_\sigma)$ for the set of active permutations of $U^+_\sigma$ with respect to $z_\zeta$ and $w$, i.e.,
\[
B_{vw}(U^+_\sigma)=\{\pi \in S_{V(U^+_\sigma)} \mid \pi(U_\sigma-\lambda^{-1}(z_\zeta))=T_\sigma-\pi(w)\}.
\]
Clearly, $U^+_\sigma$ is a sunshine graph with skeleton $z_0z_1 \ldots z_{c-1}z_0$, labelled analogously to the labelling of $U^+$, and $U_\sigma$ is a sunshine graph. Furthermore, since $d_{U^+_\sigma}(z_\zeta)=2$, it follows from Lemma \ref{sunshine lemma 1} that $T_\sigma$ is a caterpillar. Hence $U^+_\sigma$ is a supercard of the sunshine graph $U_\sigma$ and the caterpillar $T_\sigma$. We may therefore use results corresponding to Lemma \ref{supercard sunshine lemma new 1} to Corollary \ref{supercard sunshine cor 7a} by substituting $U^+_\sigma$, $U_\sigma$, $T_\sigma$ and $B_{vw}(U^+_\sigma)$ for $U^+$, $U$, $T$ and $B_{vw}(U^+)$, respectively.

We now show that $U^+_\sigma$ is also a supercard of $U$ and $T$, and thence relate the maximum saturating sets of $B_{vw}(U^+)$ and $B_{vw}(U_\sigma^+)$. In each of the following four lemmas, we define $U^+_\sigma$ to be the supercard in equation (\ref{sunshine supercard 2 defn}) for the given permutation $\sigma$, satisfying $\sigma(U-x_\xi)=T- \sigma(w)$, where $x_\xi=\sigma^{-1}(x_\nu)$. As above, $\theta$ denotes the corresponding map from $V(U^+)$ to $V(U^+_\sigma)$.
\begin{lemma}\label{supercard sunshine lem 7a}
{\rm Let $\sigma \in \widetilde{B}_{vw}(U_\sigma^+)$. Then $\theta^{-1}(U_\sigma)=U$, $\sigma\theta^{-1}(T_\sigma)=T$, so $U^+_\sigma$ is a supercard of $U$ and $T$. 

\textit{Proof} Since $\theta(w)=w$, the restriction of $\theta$ to $U$ is clearly a relabelling of the vertices of $U$ that preserves neighbourhoods. So $\theta(U)= U_\sigma$. It now follows that
\[
\sigma\theta^{-1}(U^+_\sigma-w-z_\zeta)=\sigma(U-\theta^{-1}(z_\zeta))=\sigma(U-x_\xi)=T-\sigma(w).
\]
So, since $\sigma\theta^{-1}(N_{U^+_\sigma-z_\zeta}(w))=\{\sigma\theta^{-1}(z_0)\}=\{\sigma(x_\eta)\}=N_T(\sigma(w))$, it follows from Lemma \ref{iso cards and matching nhoods lemma} that $\sigma\theta^{-1}(U^+_\sigma-z_\zeta)=T$.
\hfill\eot}
\end{lemma}
For any $\lambda \in B_{vw}(U^+)$, we define $\lambda_\sigma \in S_{V(U^+_\sigma)}$ by $\lambda_\sigma=\theta\sigma^{-1}\lambda\theta^{-1}$.
\begin{lemma}\label{supercard sunshine lem 7b1}
{\rm Let $\sigma \in \widetilde{B}_{vw}(U_\sigma^+)$ and let $\lambda \in B_{vw}(U^+)$. Then $\lambda_\sigma^{-1}(z_\zeta)=\theta(x_\mu)$ and $\lambda_\sigma \in  B_{vw}(U^+_\sigma)$.

\textit{Proof} $\lambda_\sigma^{-1}(z_\zeta)=\theta\lambda^{-1}\sigma\theta^{-1}(z_\zeta)=\theta(x_\mu)$. In addition, since $\theta^{-1}(U_\sigma)=U$, $\sigma\theta^{-1}(T_\sigma)=T$ and $\lambda(U-x_\mu)=T-\lambda(w)$, we have
\[
\lambda_\sigma(U_\sigma-\lambda_\sigma^{-1}(z_\zeta))=\lambda_\sigma(U_\sigma-\theta(x_\mu))=\theta\sigma^{-1}\lambda(U-x_\mu)=\theta\sigma^{-1}(T-\lambda(w))=T_\sigma-\theta\sigma^{-1}\lambda\theta^{-1}(w),
\]
as $\theta(w)=w$. Hence $\lambda_\sigma(U_\sigma-\lambda_\sigma^{-1}(z_\zeta))=T_\sigma-\lambda_\sigma(w)$, so $\lambda_\sigma \in  B_{vw}(U^+_\sigma)$.
\hfill\eot}
\end{lemma}
\begin{lemma}\label{supercard sunshine lem 7b}
{\rm Let $X$ be a maximum saturating set of $B_{vw}(U^+)$ and suppose that $\sigma \in \widetilde{X}$. Then the set $X_\sigma$ defined by $X_\sigma=\{\lambda_\sigma \mid \lambda \in X\}$ is a maximum saturating set of $B_{vw}(U^+_\sigma)$.

\textit{Proof} By Lemma \ref{supercard sunshine lem 7b1}, $X_\sigma \subseteq B_{vw}(U^+_\sigma)$. Moreover, since $X$ is a maximum saturating set of $B_{vw}(U^+)$ that contains $\sigma$, it is straightforward to show that $X_\sigma$ satisfies conditions (a) and (b) of Definition \ref{supercard number of common cards}. So, since $|X_\sigma|=|X|$, it follows from Theorem \ref{max sat set lemma}(c) that $X_\sigma$ is a maximum saturating set of $B_{vw}(U^+_\sigma)$.
\hfill\eot}
\end{lemma}
For the final lemma in this section, we make use of the fact that if $\theta(x_i)=z_j$ then $\theta(x_{i+k}) \in \{z_{j-k},\, z_{j+k}\}$, and $\theta^{-1}(z_{j+k}) \in \{x_{i-k},\, x_{i+k}\}$ for all $k$.
\begin{lemma}\label{supercard sunshine lem 7d}
{\rm Let $X$ be a $U^+$-optimum saturating set of $B_{vw}(U^+)$. Suppose there exists some $\sigma \in \widetilde{X} \setminus X_U$. Then $|\widetilde{X}|  \le 2 \max(\chi(U^+) ,\, \chi(U_\sigma^+))+1$.

\textit{Proof} Let $X_\sigma=\{\pi_\sigma \mid \pi \in X\}$, and let $P$ be the subset of $X_\sigma$ defined by $P =\{\pi_\sigma \mid \pi \in \widetilde{X} \setminus X_U\}$. We show that $|P| \le \chi(U_\sigma^+)+1$. As $|P|=|\widetilde{X} \setminus X_U|$, it will then follow that $|\widetilde{X}| \le \chi(U^+) + \chi(U_\sigma^+)+1 \le  2 \max(\chi(U^+) ,\, \chi(U_\sigma^+))+1$. We note that, since $X_\sigma$ is a maximum saturating set of $B_{vw}(U^+_\sigma)$ by Lemma \ref{supercard sunshine lem 7b}, we may define $\widetilde{X}_\sigma=X_\sigma \cap \widetilde{B}_{vw}(U_\sigma^+)$ as in Definition \ref{supercard sunshine defn X_tilde}.

It follows from the definition of $P$ that, given any $\pi_\sigma \in P$, there exists a corresponding $\pi \in \widetilde{X} \setminus X_U$. We next show that if $\lambda_\sigma \in P \cap \widetilde{X}_\sigma$ then $\lambda_\sigma(z_{\alpha+2}) \in \{z_{\zeta-2}),\, z_{\zeta+2}\}$, where $z_\alpha=\lambda_\sigma^{-1}(z_\zeta)$. So $P \cap \widetilde{X}_\sigma \subseteq L_\sigma$, where $L_\sigma$ is the subset of $X_\sigma$ that corresponds to the subset $L$ of $X$ in Lemma \ref{supercard sunshine lemma 7 new}. Using Lemma \ref{supercard sunshine lemma 7 new}(a) for $X_\sigma$, it will immediately follow that $|P \cap \widetilde{X}_\sigma| \le \chi(U_\sigma^+)$.

Let $\lambda_\sigma \in P \cap \widetilde{X}_\sigma$. On using Corollary \ref{supercard sunshine cor 7a} first for $\lambda$ and then for $\sigma$, we see that $\{\lambda(x_{\mu-2}), \lambda(x_{\mu+2})\}=\{\sigma(x_{\xi-2}), \sigma(x_{\xi+2})\}$. Thus $\{\sigma^{-1}\lambda(x_{\mu-2}), \sigma^{-1}\lambda(x_{\mu+2})\} = \{x_{\xi-2}, x_{\xi+2}\}$. So, since $z_\alpha=\theta(x_\mu)$ by Lemma \ref{supercard sunshine lem 7b1}, it follows that $\theta^{-1}(z_{\alpha+2}) \in \{x_{\mu-2}, x_{\mu+2}\}$, and hence $\sigma^{-1}\lambda\theta^{-1}(z_{\alpha+2}) \in \{x_{\xi-2}, x_{\xi+2}\}$. Therefore $\lambda_\sigma(z_{\alpha+2}) \in \{\theta(x_{\xi-2}),\, \theta(x_{\xi+2})\}= \{z_{\zeta-2}),\, z_{\zeta+2}\}$ as required. 

It remains to be shown that $|P \setminus \widetilde{X}_\sigma| \le 1$. So suppose that there exists $\lambda_\sigma \in P \setminus \widetilde{X}_\sigma$. Now $\lambda(w)$ is a $d$-leaf of $T$ as $\lambda \in \widetilde{X}$. So, since $\theta\sigma^{-1}(T)=T_\sigma$  by Lemma \ref{supercard sunshine lem 7a} and $\lambda_\sigma = \theta\sigma^{-1}\lambda\theta^{-1}$, it follows that $\lambda_\sigma(w)$ must be a $d$-leaf of $T_\sigma$. As $\lambda_\sigma \not\in  \widetilde{X}_\sigma$, this implies that $\lambda_\sigma(w)$ is a leaf of $T_\sigma$ but not a leaf of $U^+_\sigma$. We now show that $z_{\zeta+1}$ is the unique leaf of $T_\sigma$ that is not a leaf of $U^+_\sigma$. Since $\pi_\sigma(w)$ must be distinct for each $\pi_\sigma \in X_\sigma$, this will imply that $P \setminus \widetilde{X}_\sigma=\{\lambda_\sigma\}$, and therefore $|P \setminus \widetilde{X}_\sigma| \le 1$.

Now $\mathrm{diam}(T)=\mathrm{diam}(T_\sigma)$ as $T \cong T_\sigma$. So, by applying Lemma \ref{sunshine lemma 1} to $U^+$ and $x_\nu$, and then to $U_\sigma^+$ and $z_\zeta$, it is easy to see that $\tau_{U^+}(x_\nu)=\tau_{U_\sigma^+}(z_\zeta)$. Now, since $X$ is $U^+$-optimum and $\sigma \in \widetilde{X} \setminus X_U$, it follows from Corollary \ref{supercard sunshine cor 7a} that $\tau_{U^+}(x_\nu)=1$, and thus $\tau_{U_\sigma^+}(z_\zeta)=1$. On applying Lemma \ref{supercard sunshine lemma new 1} to $U^+_\sigma$, $T_\sigma$ and $z_\zeta$, it then follows that $z_{\zeta+1}$ is the only leaf of $T_\sigma$ that is not a leaf of $U^+_\sigma$. This completes the proof.
\hfill\eot}
\end{lemma}
\begin{thm}\label{supercard sunshine thm 7e}
{\rm Let $U^+$ be a supercard of $U$ and $T$ that is a sunshine graph that has the largest possible value of $\chi(U^+)$ over all supercards of $U$ and $T$ that are sunshine graphs. Then $b(U, \,T) \le 2\chi(U^+) + 5$.

\textit{Proof} Let $X$ be a $U^+$-optimum saturating set of $B_{vw}(U^+)$. Now, if $\widetilde{X} \subseteq X_U$ then the result holds immediately by Corollary \ref{supercard sunshine cor 4d}. So suppose that there exists $\sigma \in \widetilde{X} \setminus X_U$, and let $U^+_\sigma$ be the supercard of $U$ and $T$ as defined in equation (\ref{sunshine supercard 2 defn}). Then $|\widetilde{X}| \le 2 \max(\chi(U^+) ,\, \chi(U_\sigma^+))+1$ by Lemma \ref{supercard sunshine lem 7d}. So, since $U^+_\sigma$ is a sunshine graph, $b(U, \,T) \le 2\chi(U^+) + 5$ by Corollary \ref{supercard sunshine cor 4d}.
\hfill\eot}
\end{thm}
\vspace{-5mm}
\section{\normalsize The set $B^U_{vw}(U^+)$}\label{isomorphic cards of widehat U}\vspace{-2.5mm}
Let $U$ be a sunshine graph and $T$ be a caterpillar. {\em For the rest of this paper, we shall assume that $U^+$ is a supercard of $U$ and $T$ that satisfies the conditions of Theorem \ref{supercard sunshine thm 7e}.} In light of the bound in this theorem, we now concentrate on the set $B^U_{vw}(U^+)$. For ease of notation, we write $B_U$ instead of $B^U_{vw}(U^+)$.

Let $\lambda \in B_U$. By Lemma \ref{supercard sunshine lemma 6}, $\lambda(w)$ is a leaf of $U^+$ adjacent to a vertex of degree $d_{U^+}(x_0)$. Moreover, $\lambda$ is either a rotation of the cycle $\mathrm{skel}(U)$ and $\lambda(x_i)=x_{\nu-\mu+i}$ for each $x_i$, or $\lambda$ is a reflection of the cycle $\mathrm{skel}(U)$ and $\lambda(x_i)=x_{\nu+\mu-i}$ for each $x_i$. We may therefore partition $B_U$ into the {\em rotations} $\mathrm{Rot}(B_U)$, and the {\em reflections} $\mathrm{Ref}(B_U)$.

We make frequent use of the following well-known results concerning the rotations and reflections of a cycle.
\begin{lemma}\label{supercard sunshine iso cards lemma 1}
{\rm  Let $\lambda,\, \pi \in B_U$, and let $x_\alpha=\lambda(x_0)$ and $x_\beta=\pi(x_0)$. The following results holds for all $x_i$.
\begin{enumerate}\lettnum\vspace{-4mm}
\item If $\lambda \in \mathrm{Rot}(B_U)$ then $\lambda(x_i)=x_{\alpha+i}$ and $\lambda^{-1}(x_i)=x_{i-\alpha}$. \vspace{-2.5mm}
\item If $\lambda \in \mathrm{Ref}(B_U)$ then $\lambda(x_i)=\lambda^{-1}(x_i)=x_{\alpha-i}$. \vspace{-2.5mm}
\item If $\lambda,\,\pi \in \mathrm{Rot}(B_U)$ then $\lambda\pi(x_i)=\pi\lambda(x_i)=x_{(\alpha+\beta)+i}$. \vspace{-2.5mm}
\item If $\lambda,\, \pi \in \mathrm{Ref}(B_U)$ then $\pi\lambda(x_i)=x_{(\beta-\alpha)+i}$ and $\lambda\pi(x_i)=x_{(\alpha-\beta)+i}$.\vspace{-2.5mm}
\item If $\lambda \in \mathrm{Rot}(B_U)$ and $\pi \in \mathrm{Ref}(B_U)$ then  $\lambda\pi(x_i)=x_{(\beta+\alpha)-i}$, and $\pi\lambda(x_i)=x_{(\beta-\alpha)-i}$.\vspace{-2.5mm}
\end{enumerate}
\vspace{-8mm}
\hfill\eot}
\end{lemma}
\begin{cor}\label{supercard sunshine iso cards cor 1a}
{\rm Let $\lambda \in B_U$. Suppose there exists $\sigma \in B_U$ such that $\sigma(x_i)=\lambda^{-1}(x_i)$ for all $x_i$. Then either $\lambda$ and $\sigma$ are both in $\mathrm{Rot}(B_U)$, or they are both in $\mathrm{Ref}(B_U)$.
\hfill\eot}
\end{cor}
\begin{cor}\label{supercard sunshine iso cards cor 1b}
{\rm Let $\lambda,\, \pi$ be in $B_U$. Suppose there exists $\sigma \in B_U$ such that $\sigma(x_i)=\lambda\pi(x_i)$ for all $x_i$. Then $\sigma \in \mathrm{Rot}(B_U)$ if and only if either $\lambda$ and $\pi$ are both in $\mathrm{Rot}(B_U)$, or they are both in $\mathrm{Ref}(B_U)$; otherwise $\sigma \in \mathrm{Ref}(B_U)$.
\hfill\eot}
\end{cor}
We call any $\lambda \in \mathrm{Rot}(B_U)$ such that $\lambda(x_0)=x_0$, a {\em trivial rotation} (equivalently, $\lambda(x_i)=x_i$ for all $x_i$ by Lemma \ref{supercard sunshine iso cards lemma 1}(a)). Clearly, $\lambda^{-1}(x_\nu) = x_\nu$ for every trivial rotation $\lambda$ in $B_U$. It then immediately follows from Definition \ref{supercard number of common cards} that $1_{V(U^+)}$ is the only trivial rotation in any maximum saturating set $X$ of $B_{vw}(U^+)$.
\begin{defn}\label{X_rot X_ref}
{\rm Let $X$ be a maximum saturating set of $B_{vw}(U^+)$. We define $X_{\mathrm {Rot}}=X \cap \mathrm {Rot}(B_U)$, $X_{\mathrm {Ref}}=X \cap \mathrm {Ref}(B_U)$ and $X_{\mathrm {Aut}}=X \cap \mathrm {Aut}(U^+)$. 
\hfill\eot}
\end{defn}
\begin{lemma}\label{distinct rots and refl lemma}
{\rm Let $X$ be a maximum saturating set of $B_{vw}(U^+)$, and let $\lambda$ and $\pi$ be distinct permutations in $X$. If $\lambda$ and $\pi$ are both in $X_{\mathrm {Rot}}$ then $\lambda(x_i) \ne \pi(x_i)$ for all $x_i$. Similarly, if $\lambda$ and $\pi$ are both in $X_{\mathrm {Ref}}$ then $\lambda(x_i) \ne \pi(x_i)$ for all $x_i$.

\textit{Proof} Let $x_\alpha=\lambda(x_0)$ and $x_\beta=\pi(x_0)$. Suppose that there exists $x_k$ such that $\lambda(x_k) =\pi(x_k)$. Now, if both $\lambda$ and $\pi$ are in $\mathrm{Rot}(B_U)$ then it follows from Lemma \ref{supercard sunshine iso cards lemma 1}(a) that $\alpha=\beta$, so $\lambda^{-1}(x_\nu)=\pi^{-1}(x_\nu)$. Similarly, if both $\lambda$ and $\pi$  are in $\mathrm{Ref}(B_U)$ then it follows from Lemma \ref{supercard sunshine iso cards lemma 1}(b) that $\lambda^{-1}(x_\nu)=\pi^{-1}(x_\nu)$. Either contradicts property (b) of Definition \ref{supercard number of common cards}.
\hfill\eot}
\end{lemma}
\begin{cor}\label{distinct rots and refl cor}
{\rm Let $X$ be a maximum saturating set of $B_{vw}(U^+)$, and let $Z=\{\lambda(x_0) \mid \lambda \in B_U\}$. Then $|X_U| \le 2|Z|$.
 
\textit{Proof} For each $z \in Z$, it follows from Lemma \ref{distinct rots and refl lemma} that there is at most one $\lambda \in X_{\mathrm {Rot}}$ and at most one $\pi \in X_{\mathrm {Ref}}$ such that $\lambda(x_0)=\pi(x_0)=z$.
\hfill\eot}
\end{cor} 
We now make some further observations about $B_U$.
\begin{lemma}\label{iso cards degrees lemma}
{\rm Let $\theta \in B_U$ and let $u \in V(U^+)$. Then
\begin{enumerate}\lettnum \vspace{-4mm}
\item $d_{U^+}(u)=d_{U^+}(\theta(u))+1$ if and only if $u$ is $x_0$ and $\theta(w)$ is not adjacent to $\theta(x_0)$;\vspace{-2.5mm}
\item $d_{U^+}(u)=d_{U^+}(\theta(u))-1$ if and only if $u$ is not $x_0$ and $\theta(w)$ is adjacent to $\theta(u)$;\vspace{-2.5mm}
\item $d_{U^+}(u)=d_{U^+}(\theta(u))$ otherwise.\vspace{-2.5mm}
\end{enumerate}
\textit{Proof} Since $\theta(U)= U^+-\theta(w)$, it follows that $d_U(u)=d_{U^+-\theta(w)}(\theta(u))$ for any $u \in V(U)$. It is now easy to show that one of (a), (b) and (c) must hold as $w$ is adjacent to $x_0$ in $U^+$. 
\hfill\eot}
\end{lemma}
\begin{cor}\label{lambda(w) adj x_0 new cor}
{\rm Let $\theta \in B_U$ be such that $\theta(w)$ is adjacent to $x_0$. If $x_i$ is not $x_0$ and $\theta(x_i) \ne x_0$ then $d_{U^+}(x_i)=d_{U^+}(\theta(x_i))$.
\hfill\eot}
\end{cor}
\begin{cor}\label{iso cards degrees cor}
{\rm Let $\theta \in B_U$.
\begin{enumerate}\lettnum \vspace{-4mm}
\item $d_{U^+}(x_0)=d_{U^+}(\theta(x_0))+1$ if and only if $\theta \not\in \mathrm {Aut}(U^{+})$;\vspace{-2.5mm}
\item $d_{U^+}(x_0)=d_{U^+}(\theta(x_0))$ if and only if $\theta \in \mathrm {Aut}(U^{+})$.\vspace{-2.5mm}
\item If $\theta(x_0)=x_0$ then $\theta \in \mathrm {Aut}(U^+)$.\vspace{-2.5mm}
\end{enumerate}
\textit{Proof} By Corollary \ref{iso cards and matching nhoods aut cor}, $\theta \in  \mathrm {Aut}(U^{+})$ if and only if $\theta(w)$ is adjacent to $\theta(x_0)$. Parts (a) and (b) then follow by putting $u=x_0$ in Lemma \ref{iso cards degrees lemma} when $u$ is $x_0$. Part (c) is immediate from part (b).
\hfill\eot}
\end{cor}
It follows from Corollary \ref{iso cards degrees cor}(c) that any trivial rotation in $B_U$ is an automorphism of $U^+$.
\begin{lemma}\label{one auto one not lemma}
{\rm Let $\lambda \in  \mathrm {Aut}(U^+)$ and $\pi \in  B_U \setminus \mathrm {Aut}(U^+)$. Then  $\lambda\pi \in  B_U \setminus \mathrm {Aut}(U^+)$ and $d_{U^+}(\lambda\pi (x_0))=d_{U^+}(x_0)-1$.

\textit{Proof} Since $\pi \in B_U$ and $\lambda \in \mathrm {Aut}(U^+)$, clearly $\lambda\pi(U)=\lambda(U^+-\pi(w))=U^+-\lambda\pi(w)$, so $\lambda\pi \in B_U$. Now, if $\lambda\pi \in \mathrm {Aut}(U^+)$ then $\pi=\lambda^{-1}\lambda\pi \in \mathrm {Aut}(U^+)$ as $\mathrm {Aut}(U^+)$ is a group. Since this is impossible, $\lambda\pi \in  B_U \setminus \mathrm {Aut}(U^+)$. So  $d_{U^+}(\lambda\pi (x_0))=d_{U^+}(x_0)-1$ by Corollary \ref{iso cards degrees cor}(a).
\hfill\eot}
\end{lemma}
\begin{lemma}\label{lambda sq lemma}
{\rm Let $\lambda \in B_U\setminus \mathrm {Aut}(U^{+})$ be such that  $\lambda^2(x_0)=x_0$. Then $\lambda(w)$ is adjacent to $x_0$.

\textit{Proof} $d_{U^+}(x_0) = d_{U^+}(\lambda(x_0))+1$ by Lemma \ref{iso cards degrees cor}(a). So, since $\lambda^2(x_0)=x_0$, clearly $d_{U^+}(\lambda(x_0)) =d_{U^+}(\lambda^2(x_0))-1$. Therefore $\lambda(w)$ is adjacent to $\lambda^2(x_0)$, i.e. $x_0$, by Lemma \ref{iso cards degrees lemma}(b).
\hfill\eot}
\end{lemma}
\begin{cor}\label{lambda ref sq cor}
{\rm Let $\lambda \in \mathrm{Ref}(B_U) \setminus \mathrm {Aut}(U^+)$. Then $\lambda(w)$ is adjacent to $x_0$.

\textit{Proof} As $\lambda  \in \mathrm{Ref}(B_U)$, we have $\lambda^2(x_0)=x_0$ by Lemma \ref{supercard sunshine iso cards lemma 1}(b).
\hfill\eot}
\end{cor}
We next show that if any non-trivial rotation is an automorphism of $U^+$ then $B_U = \mathrm {Aut}(U^{+})$.
\begin{lemma}\label{sunshine supercard one rot auto all rot autos lemma}
{\rm Suppose there exists $\lambda \in \mathrm{Rot}(B_U) \cap \mathrm {Aut}(U^+)$ such that $\lambda(x_0) \ne x_0$. Then $B_U = \mathrm {Aut}(U^{+})$.

\textit{Proof} Let us assume, to the contrary, that there exists $\pi \in B_U \setminus \mathrm {Aut}(U^{+})$. Suppose first that $\pi \in \mathrm{Rot}(B_U)$. Then, by Lemma \ref{supercard sunshine iso cards lemma 1}(c) and Lemma \ref{one auto one not lemma}, $d_{U^+}(\pi\lambda(x_0)) =d_{U^+}(\lambda\pi(x_0)) =d_{U^+}(x_0)-1$. Thus $d_{U^+}(\lambda(x_0))=d_{U^+}(\pi\lambda(x_0))+1$ as $\lambda \in \mathrm {Aut}(U^{+})$. However by using Lemma \ref{iso cards degrees lemma}(a) with $\theta=\pi$ and $u=\lambda(x_0)$ yields $\lambda(x_0) = x_0$, which is impossible.

Suppose, on the other hand, that $\pi \in \mathrm{Ref}(B_U)$. Then, by Lemma \ref{one auto one not lemma} and Corollary \ref{supercard sunshine iso cards cor 1b}, $\lambda\pi \in \mathrm{Ref}(B_U) \setminus \mathrm {Aut}(U^+)$. So $\lambda\pi(w)$ is adjacent to $x_0$ by Corollary \ref{lambda ref sq cor}. By Lemma \ref{supercard sunshine iso cards lemma 1}(b), $\lambda\pi(\pi(x_0))=\lambda(x_0)$. So, since $\pi(x_0) \ne x_0$ by Corollary \ref{iso cards degrees cor}(c) and $\lambda(x_0) \ne x_0$ by assumption, it follows from Corollary \ref{lambda(w) adj x_0 new cor} with $\theta=\lambda\pi$ and $x_i=\pi(x_0)$ that $d_{U^+}(\pi(x_0))=d_{U^+}(\lambda\pi(\pi(x_0)))=d_{U^+}(\lambda(x_0))$. However, Corollary \ref{iso cards degrees cor}(a) and (b) applied to $\pi$ and $\lambda$ respectively, implies that $d_{U^+}(\pi(x_0)) \ne d_{U^+}(\lambda(x_0))$, contradicting the previous statement. Hence $B_U= \mathrm {Aut}(U^{+})$.
\hfill\eot}
\end{lemma}
\begin{cor}\label{sunshine supercard two distinct auto refs all rot autos cor 1}
{\rm Let $X$ be a $U^+$-optimum saturating set of $B_{vw}(U^+)$.
\begin{enumerate}\lettnum\vspace{-4mm}
\item If $B_U \ne \mathrm{Aut}(U^+)$ then $X_{\mathrm{Aut}} \subseteq X_{\mathrm{Ref}} \cup \{1_{V(U^+)} \}$.\vspace{-2.5mm}
\item If $|X_{\mathrm{Aut}}| \ge 3$ then there exists a non-trivial rotation in $\mathrm{Aut}(U^+)$ and $B_U=\mathrm{Aut}(U^+)$.\vspace{-2.5mm}
\item If $b(U,\, T) \ge 10$ and $B_U=\mathrm{Aut}(U^+)$ then there exists a non-trivial rotation in $\mathrm{Aut}(U^+)$.\vspace{-2.5mm}
\end{enumerate}
\textit{Proof} (a) $1_{V(U^+)}$ is the only trivial rotation in $X$. So, if there exists another rotation in $X_{\mathrm{Aut}}$ then $B_U = \mathrm{Aut}(U^+)$ by Lemma \ref{sunshine supercard one rot auto all rot autos lemma}.

(b) Suppose that $\lambda$ and $\pi$ are distinct permutations in $X_{\mathrm{Aut}} \setminus \{1_{V(U^+)}\}$. By Lemma \ref{sunshine supercard one rot auto all rot autos lemma}, we may assume that both $\lambda$ and $\pi$ are reflections as $1_{V(U^+)}$ is the only trivial rotation in $X$. Now, since $\mathrm{Aut}(U^+)$ is a group, it follows from Corollary \ref{supercard sunshine iso cards cor 1b} that $\lambda\pi$ is a rotation in $\mathrm {Aut}(U^+)$. Moreover, since $\lambda(x_0) \ne \pi(x_0)$ by Lemma \ref{distinct rots and refl lemma}, $\lambda\pi(x_0) \ne x_0$ by Lemma \ref{supercard sunshine iso cards lemma 1}(d). Therefore $B_U= \mathrm {Aut}(U^+)$ by Lemma \ref{sunshine supercard one rot auto all rot autos lemma}.

(c) $|X_{\mathrm{Aut}}| = \chi(U^+)$ as $X$ is $U^+$-optimum and $B_U=\mathrm{Aut}(U^+)$. The result now follows from Theorem \ref{supercard sunshine thm 7e} and part (b).
\hfill\eot}
\end{cor}
We now consider those rotations that are not automorphisms of $U^+$.
\begin{lemma}\label{rot Aut(U) lemma 1}
{\rm Suppose there exists $\lambda,\, \pi \in \mathrm{Rot}(B_U) \setminus \mathrm {Aut}(U^+)$ such that $\lambda(x_0) \ne \pi(x_0)$. Then both $\lambda(w)$ and $\pi(w)$ are adjacent to $x_0$.

\textit{Proof} Let $x_\alpha=\lambda(x_0)$ and $x_\beta=\pi(x_0)$, so $\lambda(x_i)=x_{\alpha+i}$ and $\pi(x_i)=x_{\beta+i}$ by Lemma \ref{supercard sunshine iso cards lemma 1}(a). Now, by Corollary \ref{iso cards degrees cor}(c), $\alpha \ne 0$ as $\lambda \not \in \mathrm {Aut}(U^{+})$. So, by Lemma \ref{iso cards degrees lemma}, $d_{U^{+}}(x_{-\alpha}) \le d_{U^{+}}(\pi(x_{-\alpha})) = d_{U^{+}}(x_{\beta-\alpha})$.

Suppose that $\lambda(w)$ is not adjacent to $x_0$, i.e. $\lambda(x_{-\alpha})$. Then, since $\alpha \ne 0$, it follows from Lemma \ref{iso cards degrees lemma} that $d_{U^{+}}(x_{-\alpha}) =  d_{U^{+}}(\lambda(x_{-\alpha}))=d_{U^{+}}(x_0)$. So $d_{U^{+}}(x_{\beta-\alpha}) \ge d_{U^{+}}(x_0)$. Now, by Corollary \ref{iso cards degrees cor}(a), we have that $ d_{U^{+}}(x_0) =  d_{U^{+}}(\pi(x_0))+1$ as $\pi \not\in \mathrm {Aut}(U^{+})$. Thus
\[
d_{U^{+}}(x_{\beta-\alpha})  \ge d_{U^{+}}(x_0) >  d_{U^{+}}(\pi(x_0))= d_{U^{+}}(x_\beta)=d_{U^{+}}(\lambda(x_{\beta-\alpha})).
\]
Hence, $x_{\beta-\alpha}$ must be $x_0$ by Lemma \ref{iso cards degrees lemma}(a), and therefore, $\beta=\alpha$, which is impossible. So $\lambda(w)$ must be adjacent to $x_0$. By symmetry, clearly $\pi(w)$ is also adjacent to $x_0$.
\hfill\eot}
\end{lemma}
\begin{cor}\label{rot Aut(U) cor 3a}
{\rm Let $X$ be a maximum saturating set of $B_{vw}(U^+)$. If $|X_{\mathrm{Rot}} \setminus X_{\mathrm{Aut}}| \ge 2$ then $\lambda(w)$ is adjacent to $x_0$ for all $\lambda \in X_{\mathrm{Rot}} \setminus X_{\mathrm{Aut}}$.

\textit{Proof} Suppose that there exist distinct rotations $\lambda$ and $\pi$ in $X_{\mathrm{Rot}} \setminus X_{\mathrm{Aut}}$. Then $\lambda(x_0) \ne \pi(x_0)$ by Lemma \ref{distinct rots and refl lemma}, so the result follows immediately from Lemma \ref{rot Aut(U) lemma 1}.
\hfill\eot}
\end{cor}
The following lemma shows that there exists a correspondence between those permutations $\lambda$ for which $\lambda(w)$ is adjacent to $x_0$, and a subset of $\mathrm{Aut}(U)$.
\begin{lemma}\label{B_U Aut(U) lemma 1}
{\rm \begin{enumerate}\lettnum\vspace{-2.5mm}
\item For each $\lambda \in B_U$ for which the leaf $\lambda(w)$ of $U^+$ is adjacent to $x_0$, there exists $\lambda^\ast \in \mathrm {Aut}(U)$ such that $\lambda^\ast(x_i)=\lambda(x_i)$ for all $x_i$.\vspace{-2.5mm}
\item For each $\lambda^\ast \in \mathrm {Aut}(U)$, there exists $\lambda \in B_U$ such $\lambda(w)=w$ and $\lambda(x_i)=\lambda^\ast(x_i)$ for all $x_i$.
\end{enumerate}
\textit{Proof}
(a) Let $\lambda \in B_U$ be such that $\lambda(w)$ is adjacent to $x_0$. Let $\phi$ be the transposition that swaps the leaves $w$ and $\lambda(w)$, and let $\lambda^\ast=\phi\lambda$. Then $\lambda^\ast(U)=\phi\lambda(U)=\phi(U^+-\lambda(w))=U$, as $\phi \in \mathrm {Aut}(U^+)$.

(b) Let $\lambda^\ast \in \mathrm {Aut}(U)$. Then the permutation $\lambda$ defined by $\lambda(w)=w$ and $\lambda(u)=\lambda^\ast(u)$ for all $u \in V(U)$ clearly has the required properties.
\hfill\eot}
\end{lemma}
\begin{lemma}\label{B_U Aut(U) lemma 2}
{\rm Suppose there exists $\lambda,\, \pi \in \mathrm{Ref}(B_U) \setminus \mathrm {Aut}(U^+)$ such that $\lambda(x_0) \ne \pi(x_0)$. Then there exists $\sigma \in \mathrm{Rot}(B_U) \setminus \mathrm {Aut}(U^+)$ such that $\sigma(w)=w$ and $\sigma(x_i)=\lambda\pi(x_i)$ for all $x_i$.

\textit{Proof} Since $\lambda(w)$ and $\pi(w)$ are both adjacent to $x_0$ by Corollary \ref{lambda ref sq cor}, it follows from Lemma \ref{B_U Aut(U) lemma 1}(a) that there exist $\lambda^\ast$ and $\pi^\ast$ in $\mathrm{Aut}(U)$ such that $\lambda^\ast(x_i)=\lambda(x_i)$ and $\pi^\ast(x_i)=\pi(x_i)$ for all $x_i$. Since $\mathrm{Aut}(U)$ is a group, $\lambda^\ast\pi^\ast \in \mathrm{Aut}(U)$. So, by Lemma \ref{B_U Aut(U) lemma 1}(b), there exists $\sigma \in B_U$ such that $\sigma(w)=w$ and $\sigma(x_i)=\lambda\pi(x_i)$ for all $x_i$. 

Now, by Corollary \ref{supercard sunshine iso cards cor 1b}, $\sigma \in \mathrm{Rot}(B_U)$. Moreover, since $\lambda(x_0) \ne \pi(x_0)$, it follows from Lemma \ref{supercard sunshine iso cards lemma 1}(d) that $\sigma(x_0) \ne x_0$. As $B_U \ne \mathrm{Aut}(U^+)$, it immediately follows from Lemma \ref{sunshine supercard one rot auto all rot autos lemma} that $\sigma \not\in \mathrm{Aut}(U^+)$.
\hfill\eot}
\end{lemma}
\begin{lemma}\label{rot Aut(U) lemma 3}
{\rm Suppose there exists $\lambda,\, \pi \in \mathrm{Ref}(B_U) \setminus \mathrm {Aut}(U^+)$ such that $\lambda(x_0) \ne \pi(x_0)$. Then $\psi(w)$ is adjacent to $x_0$ for all $\psi \in B_U$.

\textit{Proof} Let $\sigma$ be the rotation from Lemma \ref{B_U Aut(U) lemma 2} and let $\psi \in B_U$. If $\psi(x_0)=x_0$ then $\psi \in \mathrm{Aut}(U^+)$ by Corollary \ref{iso cards degrees cor}(c), so $\psi(w)$ is adjacent to $x_0$. We therefore assume that $\psi(x_0) \ne x_0$. We first show that $\psi \not\in\mathrm{Aut}(U^+)$. 

Suppose that $\psi \in \mathrm{Aut}(U^+)$. Since $\psi(x_0) \ne x_0$ and $B_U \ne \mathrm{Aut}(U^+)$, it follows from Lemma \ref{sunshine supercard one rot auto all rot autos lemma} that $\psi$ must be a reflection.  By Lemma \ref{one auto one not lemma} and Corollary \ref{supercard sunshine iso cards cor 1b}, we see that $\psi\sigma \in \mathrm{Ref}(B_U) \setminus \mathrm{Aut}(U^+)$. So $\psi\sigma(w)$ is adjacent to $x_0$ by Corollary \ref{lambda ref sq cor}. However, $\psi\sigma(w)=\psi(w)$ as $\sigma(w)=w$, so $\psi\sigma(w)$ must be adjacent to $\psi(x_0)$ as $\psi \in \mathrm{Aut}(U^+)$. Since $\psi(x_0) \ne x_0$ this is impossible. This contradiction shows that $\psi \not\in \mathrm{Aut}(U^+)$.

It remains to be shown that $\psi(w)$ is adjacent to $x_0$. If $\psi \in \mathrm{Ref}(B_U)$ then the result follows from Corollary \ref{lambda ref sq cor}. So suppose that $\psi \in \mathrm{Rot}(B_U) \setminus \mathrm{Aut}(U^+)$. Since the result follows from Lemma \ref{rot Aut(U) lemma 1} when $\sigma(x_0) \ne \psi(x_0)$, we may assume that $\sigma(x_0)=\psi(x_0)$, thus $\psi^{-1}(x_0)=\sigma^{-1}(x_0)$ by Lemma \ref{supercard sunshine iso cards lemma 1}(a). Hence $\sigma^{-1}(x_0) \ne x_0$. As $\sigma(w)=w$, it now follows from Lemma \ref{iso cards degrees lemma}(b) with $\theta=\sigma$ and $u=\sigma^{-1}(x_0)$, that $d_{U^+}(\sigma^{-1}(x_0))=d_{U^+}(\sigma\sigma^{-1}(x_0))-1$. So $d_{U^+}(\sigma^{-1}(x_0))=d_{U^+}(x_0)-1$, and thus, $d_{U^+}(\psi^{-1}(x_0))=d_{U^+}(x_0)-1$. The result now follows by again applying Lemma \ref{iso cards degrees lemma}(b) to $\psi$ and $\psi^{-1}(x_0)$.
\hfill\eot}
\end{lemma}
\begin{cor}\label{ref Aut(U) cor 3a}
{\rm Let $X$ be a maximum saturating set of $B_{vw}(U^+)$. If $|X_{\mathrm{Ref}} \setminus X_{\mathrm{Aut}}| \ge 2$ then $\lambda(w)$ is adjacent to $x_0$ for all $\lambda \in B_U$.

\textit{Proof} Suppose that $\lambda$ and $\pi$ are two distinct reflections in $X_{\mathrm{Ref}} \setminus X_{\mathrm{Aut}}$. Then $\lambda(x_0) \ne \pi(x_0)$ by Lemma \ref{distinct rots and refl lemma}, so the result follows immediately from Lemma \ref{rot Aut(U) lemma 3}.
\hfill\eot}
\end{cor}
We note that Corollary \ref{ref Aut(U) cor 3a} is a stronger result than the analogous result for rotations in Corollary \ref{rot Aut(U) cor 3a}, since the conclusions of Corollary \ref{ref Aut(U) cor 3a} applies to all permuations in $B_U$.
\begin{lemma}\label{X_U not autos lemma}
{\rm Suppose that $B_U \ne \mathrm {Aut}(U^+)$ and $\chi(U^+) \ge 5$. Let $X$ be a $U^+$-optimum saturating set of $B_{vw}(U^+)$. Then
\begin{enumerate}\lettnum\vspace{-4mm}
\item for each $\lambda \in X_U$, there exists a distinct leaf of $U^+$, namely $\lambda(w)$, adjacent to $x_0$;\vspace{-2.5mm}
\item $d_{U^+}(x_0) \ge \chi(U^+)+2$;\vspace{-2.5mm}
\item for each $\lambda \in X_U$, there exists $\lambda^\ast \in \mathrm{Aut}(U)$ such that $\lambda^\ast(x_0)=\lambda(x_0)$.
\vspace{-2.5mm}
\end{enumerate}
\textit{Proof} Part (b) follows immediately from part (a). In addition, part (c) follows from part (a) by Lemma \ref{B_U Aut(U) lemma 1}(a). Now, by property (b) of Definition \ref{supercard number of common cards} and Lemma \ref{supercard sunshine lemma 6}, $\lambda(w)$ is a distinct leaf of $U^+$ for each distinct $\lambda  \in X_U$. To complete the proof, we now show that each such $\lambda(w)$ is adjacent to $x_0$.

By Corollary \ref{ref Aut(U) cor 3a}, there is nothing to prove when $|X_{\mathrm {Ref}} \setminus X_{\mathrm{Aut}}| \ge 2$; so we assume that $|X_{\mathrm {Ref}} \setminus X_{\mathrm{Aut}}| \le 1$. Furthermore, since $B_U \ne \mathrm {Aut}(U^{+})$, it follows from Corollary \ref{sunshine supercard two distinct auto refs all rot autos cor 1} that $X_{\mathrm{Aut}} \subseteq X_{\mathrm{Ref}} \cup \{1_{V(U^+)} \}$ and $|X_{\mathrm{Aut}}| \le 2$. Since $|X_U| \ge 5$, this implies that $|X_{\mathrm {Rot}} \setminus X_{\mathrm{Aut}}| \ge 2$. It now follows from Corollaries \ref{rot Aut(U) cor 3a} and \ref{lambda ref sq cor} that $\lambda(w)$ is adjacent to $x_0$ for each $\lambda \in X_U \setminus X_{\mathrm{Aut}}$. As $1_{V(U^+)}(w)$ is clearly also adjacent to $x_0$, this concludes the proof in the case that $X_{\mathrm {Ref}} \cap X_{\mathrm{Aut}} = \emptyset$.

Suppose then there exists $\theta \in X_{\mathrm {Ref}} \cap X_{\mathrm{Aut}}$, and let $\phi$ and $\eta$ be permutations in $X_{\mathrm {Rot}} \setminus X_{\mathrm{Aut}}$. By Lemma \ref{one auto one not lemma} and Corollary \ref{supercard sunshine iso cards cor 1b}, both $\theta\phi$ and $\theta\eta$ are in $\mathrm{Ref}(B_U) \setminus \mathrm {Aut}(U^+)$. Furthemore, $\theta\phi(x_0) \ne \theta\eta(x_0)$ as $\phi(x_0) \ne \eta(x_0)$ by Corollary \ref{distinct rots and refl lemma}. By applying Lemma \ref{rot Aut(U) lemma 3} to $\theta\phi$ and $\theta\eta$, this implies that $\psi(w)$ is adjacent to $x_0$ for all $\psi \in B_U$.
\hfill\eot}
\end{lemma}
\begin{thm}\label{X_U not autos thm final}
{\rm Suppose that $B_U \ne \mathrm {Aut}(U^+)$ and $n \ge 12$. Then $b(U,\, T) \le 2\left\lfloor \sqrt{2n+1}\right\rfloor+3$.

\textit{Proof} We recall that $d_1(U)$ and $d_2(U)$ are the number of vertices of $U$ of degrees $1$ and $2$, respectively. Let $d_q(U)=n-d_1(U)-d_2(U)$, i.e., the number of vertices of $U$ of degree $3$ or more.

Now, if $\chi(U^+) \le 4$ then $b(U,\, T) \le 13$ by Theorem \ref{supercard sunshine thm 7e}, thus the bound holds. We may therefore assume that $\chi(U^+) \ge 5$, so the conclusions of Lemma \ref{X_U not autos lemma} hold. It follows from part(b) of that lemma  that $d_U(x_0) = d_{U^+}(x_0)-1 \ge \chi(U^+)+1 \ge 6$.

Let $X$ be a $U^+$-optimum saturating of $U^+$, and let $Z=\{\lambda(x_0) \mid \lambda \in X_U\}$. By Lemma \ref{X_U not autos lemma}(c), there exists a subset $\Lambda^\ast$ of $\mathrm{Aut}(U)$ such that $\{\lambda^\ast(x_0) \mid \l\lambda^\ast \in \Lambda^\ast\}=Z$. This implies that $U$ contains at least $|Z|$ vertices of degree $d_U(x_0)$, and hence $d_q(U) \ge |Z|$. In addition, since each vertex $x_i$ in $Z$ is adjacent to precisely $d_U(x_0)-2$ leaves, it follows that $d_1(U) \ge |Z|(d_U(x_0)-2)$. Therefore $d_1(U) \ge |Z|(\chi(U^+)-1)$.

Now, by property (b) of Definition \ref{supercard number of common cards} and Lemma \ref{supercard sunshine lemma 6}, for each distinct $\lambda$ in $X_U$ there exists a distinct vertex of $U$ of degree $2$, namely $\lambda^{-1}(x_\nu)$. So $d_2(U) \ge |X_U| = \chi(U^+) $. 
Since each vertex in $Z$ has degree at least $6$ in $U$, it therefore follows that
\[
n=|V(U)| = d_1(U)+d_2(U) +d_q(U)  \ge |Z|(\chi(U^+)-1)+\chi(U^+) + |Z| = (|Z|+1)\chi(U^+).
\]
So, since $\chi(U^+) =|X_U| \le 2|Z|$ by Corollary \ref{distinct rots and refl cor}, we have $2n \ge \chi(U^+)(\chi(U^+)+2)$. Solving for $\chi(U^+)$, yields $\chi(U^+) \le \sqrt{2n+1}-1$. Therefore $b(U,\, T) \le 2\left\lfloor \sqrt{2n+1}-1\right\rfloor+5$ by Theorem \ref{supercard sunshine thm 7e}, yielding the bound.
\hfill\eot}
\end{thm}
\vspace{-6mm}
\section{\normalsize The case when $B^U_{vw}(U^+) = \mathrm{Aut}(U^+)$}\label{B(U^+) = U^+}\vspace{-2.5mm}
{\em For the whole of this section, we assume that that there exists some non-trivial rotation in $\mathrm{Aut}(U^+)$}. So $B^U_{vw}(U^+) = \mathrm {Aut}(U^+)$ by Lemma \ref{sunshine supercard one rot auto all rot autos lemma}, and thus $X_U=X_{\mathrm {Aut}} =X_{\mathrm{Rot}} \cup X_{\mathrm{Ref}}$, for any maximum saturating set $X$. We note that, if $B^U_{vw}(U^+) = \mathrm {Aut}(U^+)$ but there does not exist such a rotation, then $b(U,\, T) \le 9$ by Corollary \ref{sunshine supercard two distinct auto refs all rot autos cor 1}(c). This motivates our assumption.

For clarity, we write $\mathrm {Rot}(\mathrm{Aut}(U^+))$ and $\mathrm {Ref}(\mathrm{Aut}(U^+))$ instead of $\mathrm {Rot}(B_U)$ and $\mathrm {Ref}(B_U)$, respectively. Since $\mathrm {Aut}(U^+)$ is a group, we can simplify Corollaries \ref{supercard sunshine iso cards cor 1a} and \ref{supercard sunshine iso cards cor 1b} as follows.
\begin{cor}\label{Aut rot ref cor 1a}
{\rm Let $\lambda \in \mathrm {Aut}(U^+)$. Then either $\lambda$ and $\lambda^{-1}$ are both in $\mathrm{Rot}(\mathrm {Aut}(U^+))$, or they are both in $\mathrm{Ref}(\mathrm {Aut}(U^+))$.
\hfill\eot}
\end{cor}
\begin{cor}\label{Aut rot ref cor 1b}
{\rm Let $\lambda,\,\pi \in \mathrm {Aut}(U^+)$. Then $\lambda\pi \in \mathrm{Rot}(\mathrm {Aut}(U^+))$ if and only if either $\lambda$ and $\pi$ are both in $\mathrm{Rot}(\mathrm {Aut}(U^+))$, or they are both in $\mathrm{Ref}(\mathrm {Aut}(U^+))$; otherwise $\lambda\pi \in \mathrm{Ref}(\mathrm {Aut}(U^+))$.

\vspace{-2.5mm}
\hfill\eot}
\end{cor}
We recall that $x_i^j$ is the $j^{th}$ distinct leaf adjacent to $x_i$. The following results two are easy to prove.
\begin{lemma}\label{Aut lemma 1}
{\rm There exists $\lambda \in \mathrm{Rot}(\mathrm{Aut}(U^+))$, where $\lambda(x_i)=x_{\alpha+i}$ and $\lambda(x_i^j)=x_{\alpha+i}^j$ for all $i$ and $j$, if and only if $d_{U^+}(x_i)=d_{U^+}(x_{\alpha+i})$ for all $i$. 
\hfill\eot}
\end{lemma}
\begin{lemma}\label{Aut lemma 1a}
{\rm There exists $\lambda \in \mathrm{Ref}(\mathrm{Aut}(U^+))$, where $\pi(x_i)=x_{\beta-i}$ and $\lambda(x_i^j)=x_{\beta-i}^j$ for all $i$ and $j$, if and only if $d_{U^+}(x_i)=d_{U^+}(x_{\beta-i})$ for all $i$.
\hfill\eot}
\end{lemma}
If $Z \subseteq \mathrm{Aut}(U^+)$, we define $Z(u) =\{\lambda(u) \mid \lambda \in Z\}$. We note that if $Z$ is a subgroup of $\mathrm{Aut}(U^+)$ then $Z(u)$ is the {\em orbit} of $u$ under the group action of $Z$ on $V(U^+)$.
\begin{lemma}\label{X_U cyclic auto subgroup lemma 0}
{\rm Let $A=\{\beta \mid x_\beta \in \mathrm {Rot}(\mathrm{Aut}(U^+))(x_0)\}$, and let $\delta$ be the smallest positive element in $A$. Then
\begin{enumerate}\lettnum\vspace{-4mm}
\item $\delta$ divides every element of $A$;\vspace{-2.5mm}
\item $2 \le \delta \le \f{c}{2}$;\vspace{-2.5mm}
\item there exists $\phi \in \mathrm {Rot}(\mathrm{Aut}(U^+))$ such that $\phi^{\f{c}{\delta}}=1_{V(U^+)}$, $\phi(x_i)=x_{\delta+i}$ and $\phi(x_i^j)=x_{\delta+i}^j$ for all $i,\, j$.
\end{enumerate}\lettnum
\textit{Proof} By assumption, there exists some non-trivial rotation in $\mathrm{Aut}(U^+)$; so $0 < \delta < c$. Let $\sigma  \in \mathrm {Rot}(\mathrm{Aut}(U^+))$ be such that $\sigma(x_0)=x_\delta$. 

(a)  Suppose that there exists $\beta \in A$ such that $\delta$ does not divide $\beta$, and let $\pi(x_0)=x_\beta$. Let $d$ be the highest common factor of $\delta$ and $\beta$. By the Euclidean algorithm, there exist integers $a$ and $b$ such that $b\delta+a\beta=d$. Now, on using  Lemma \ref{supercard sunshine iso cards lemma 1}, and Corollaries \ref{Aut rot ref cor 1a} and \ref{Aut rot ref cor 1b} repeatedly, it is easy to see that $\sigma^b\pi^a \in \mathrm{Rot}(\mathrm{Aut}(U^+))$ and $\sigma^b\pi^a(x_i)=x_{b\delta+a\beta+i}=x_{d+i}$ for all $x_i$. Therefore $d_{U^+}(x_i)=d_{U^+}(x_{d+i})$ for all $x_i$, as $\sigma^b\pi^a$ is in $\mathrm{Aut}(U^+)$. It now follows from Lemma \ref{Aut lemma 1} that there exists $\lambda \in \mathrm {Rot}(\mathrm{Aut}(U^+))$ such that $\lambda(x_0)=x_d$, and thus $d \in A$. Since $0 < d < \delta$, this contradicts the minimality of $\delta$.

(b) $\delta$ divides $c$ as $c \in A$, thus $\delta \le \f{c}{2}$. So suppose that $\delta=1$. Then, on using Corollary \ref{Aut rot ref cor 1b} and Lemma \ref{supercard sunshine iso cards lemma 1}(c) repeatedly, we see that $\sigma^\nu(x_0)=x_\nu$. This is impossible as $d_{U^+}(x_0) \ge 3$ and $d_{U^+}(x_\nu) =2$. Therefore $2 \le \delta \le \f{c}{2}$.

(c) $d_{U^+}(x_i)=d_{U^+}(x_{\delta+i})$ for all $x_i$, as $\sigma  \in \mathrm{Rot}(\mathrm{Aut}(U^+))$. Hence by Lemma \ref{Aut lemma 1}, there exists $\phi \in \mathrm {Rot}(\mathrm{Aut}(U^+))$ such that $\phi(x_i)=x_{\delta+i}$ and $\phi(x_i^j)=x_{\delta+i}^j$ for all $i,\, j$. As $\delta$ divides $c$, it is easy to see that $\phi^{\f{c}{\delta}}(x_i)=x_i$ and $\phi^{\f{c}{\delta}}(x_i^j)=x_i^j$ for all $i,\, j$. So $\phi^{\f{c}{\delta}}=1_{V(U^+)}$.
\hfill\eot}
\end{lemma}
By Corollary \ref{Aut rot ref cor 1b}, the composition of two rotations is also a rotation. We may therefore make the following definition.
\begin{defn}\lettnum\label{Phi and delta}
{\rm Let $\delta$ be the positive integer and $\phi$ be the rotation from Lemma \ref{X_U cyclic auto subgroup lemma 0}. We define $\Phi$ to be the cyclic subgroup of $\mathrm{Aut}(U^+)$ of order $\f{c}{\delta}$ generated by $\phi$, i.e., $\Phi=\{\phi^j \mid  0 \le j < \f{c}{\delta}\}$.

\vspace{-2.5mm}
\hfill\eot}
\end{defn}
{\em For the rest of this section we assume that on that $\delta$, $\phi$ and $\Phi$ are as in Definition \ref{Phi and delta}}. Since the cycle length $c=\delta|\Phi|$, and the number of leaves of $U^+$ is $d_1(U^+)$, we see that
\begin{equation}
n+1=\delta|\Phi|+d_1(U^+). \label{X_U cyclic auto subgroup lem 7 eqn 1}
\end{equation}
For any $\pi \in \mathrm{Ref}(\mathrm{Aut}(U^+))$, we denote the right and left cosets of $\Phi$ with respect to $\pi$ by $\Phi\pi$ and $\pi\Phi$, respectively. It follows from Corollary \ref{Aut rot ref cor 1b} that $\Phi\pi \subseteq \mathrm{Ref}(\mathrm{Aut}(U^+))$ and $\pi\Phi \subseteq \mathrm{Ref}(\mathrm{Aut}(U^+))$.

The orbit $\Phi(u)$ of a vertex $u$ of $U^+$ under $\Phi$ is the set $\{\phi^j(u) \mid  0 \le j < \f{c}{\delta}\}$. It follows from a well-known result from Group Theory that, for any two vertices $u$ and $t$ of $U^+$, either $\Phi(u)=\Phi(t)$ or $\Phi(u) \cap \Phi(t)=\emptyset$. Thus $t$ is in $\Phi(u)$ if and only if $\Phi(u)=\Phi(t)$. Now, for every vertex $u$ of $U^+$, clearly $1_{V(U^+)}$ is the only element of $\Phi$ that fixes $u$. It therefore follows from the Orbit-Stabiliser and Lagrange Theorems \cite{GPBK} that $|\Phi(u)|=|\Phi|=|\Phi\pi|=\f{c}{\delta}$.
\begin{lemma}\label{X_U cyclic auto subgroup lemma 4}
{\rm Let $x_i \in V(\mathrm{skel}(U^+))$. Then
\begin{enumerate}\lettnum\vspace{-4mm}
\item $\mathrm{Rot}(\mathrm{Aut}(U^+))(x_i) = \Phi(x_i)$;\vspace{-2.5mm}
\item if $\lambda \in \mathrm{Rot}(\mathrm{Aut}(U^+))$ then $\lambda^{-1}(x_i) \in \Phi(x_i)$;\vspace{-2.5mm}
\item if $\pi \in \mathrm{Ref}(\mathrm{Aut}(U^+))$ then $\Phi\pi(x_i)=\pi\Phi(x_i)$;\vspace{-2.5mm}
\item if $\mathrm{Ref}(\mathrm{Aut}(U^+))(x_i) \cap \Phi(x_i) \ne \emptyset$ then $\mathrm{Ref}(\mathrm{Aut}(U^+))(x_i) \subseteq \Phi(x_i)$.\vspace{-2.5mm}
\end{enumerate}
\textit{Proof} (a) Clearly $\Phi(x_i) \subseteq \mathrm{Rot}(\mathrm{Aut}(U^+))(x_i)$. So let $\lambda \in \mathrm{Rot}(\mathrm{Aut}(U^+))$ and let $x_\alpha=\lambda(x_0)$. Now, by Lemma \ref{X_U cyclic auto subgroup lemma 0}, there exists some positive integer $k$ such that $\alpha=k\delta$. Moreover, on using Corollary \ref{Aut rot ref cor 1b} and Lemma \ref{supercard sunshine iso cards lemma 1}, we see that $\lambda(x_i)=x_{k\delta+i}=\phi^k(x_i)$. Therefore $\lambda(x_i)  \in \Phi(x_i)$, and thus 
$\mathrm{Rot}(\mathrm{Aut}(U^+))(x_i) \subseteq \Phi(x_i)$.

(b) This follows immediately from part (a) as $\lambda^{-1} \in \mathrm{Rot}(\mathrm{Aut}(U^+))$ for all $\lambda \in  \mathrm{Rot}(\mathrm{Aut}(U^+))$ by Corollary \ref{Aut rot ref cor 1a}.

(c) Let $x_\beta=\pi(x_0)$, and let $0 \le j <\f{c}{\delta}$. Then, on using Corollary \ref{Aut rot ref cor 1b} and Lemma \ref{supercard sunshine iso cards lemma 1}, we see that $\phi^j\pi(x_i)=x_{j\delta+\beta-i}=\pi\phi^{-j}(x_i)$. So $\phi^j\pi(x_i) \in \pi\Phi(x_i)$ by part (b), and therefore $\Phi\pi(x_i) \subseteq \pi\Phi(x_i)$. It similarly follows that $\pi\Phi(x_i) \subseteq \Phi\pi(x_i)$.

(d) Suppose there exists $\pi \in \mathrm{Ref}(\mathrm{Aut}(U^+))$ such that $\pi(x_i) \in \Phi(x_i)$. Then $\Phi\pi(x_i)=\Phi(x_i)$, so $\pi\Phi(x_i)=\Phi(x_i)$ by part (c). Now let $\lambda \in \mathrm{Ref}(\mathrm{Aut}(U^+))$. Then $\pi\lambda(x_i) \in \Phi(x_i)$ by Corollary \ref{Aut rot ref cor 1b} and part (a). So $\pi\lambda(x_i) \in \pi\Phi(x_i)$, and therefore $\lambda(x_i) \in \Phi(x_i)$, yielding the result.
\hfill\eot}
\end{lemma}
We now show that we can always find some maximum saturating set $X$ such that $X_{\mathrm{Aut}}$ is isomorphic to a subgroup of $\mathrm{Aut}(U^+)$.
\begin{thm}\label{X_U cyclic auto subgroup thm 6}
{\rm Suppose that $\mathrm{Aut}(U^+)$ contains a non-trivial rotation. We consider the following two (not necessarily disjoint) cases:
\begin{enumerate}\vspace{-4mm}
\item[(i)] $\mathrm {Ref}(\mathrm {Aut}(U^+))(x_\nu) \subseteq \Phi(x_\nu)$;\vspace{-2.5mm}
\item[(ii)] $d_{U^+}(x_0) =3$ and $\mathrm {Ref}(\mathrm {Aut}(U^+))(x_0) \subseteq \Phi(x_0)$.\vspace{-2.5mm}
\end{enumerate}\vspace{-1mm}
Then
\begin{enumerate}\lettnum\vspace{-2.5mm}
\item if either case (i) or case (ii) holds, there exists a $U^+$-optimum saturating set $X$ of $B_{vw}(U^+)$ such that $X_{\mathrm{Aut}}= \Phi \cong  C(\f{c}{\delta})$, the cyclic group of order $\f{c}{\delta}$.\vspace{-2.5mm}
\item if neither case holds, there exists $\pi \in \mathrm {Ref}(\mathrm{Aut}(U^+))$ and a $U^+$-optimum saturating set $X$ of $B_{vw}(U^+)$ such that $X_{\mathrm{Aut}} =\Phi  \cup \Phi\pi \cong D(\f{2c}{\delta})$, the dihedral group of order $\f{2c}{\delta}$.\vspace{-2.5mm}
\end{enumerate}
\textit{Proof} Since $\mathrm{Aut}(U^+)$ contains a non-trivial rotation, we may define $\Phi$ as in Definition \ref{Phi and delta}. Now, by parts (a) and (b) of Lemma \ref{X_U cyclic auto subgroup lemma 4}, $\lambda(x_0) \in \Phi(x_0)$ and $\lambda^{-1}(x_\nu) \in \Phi(x_\nu)$ for all $\lambda \in \mathrm{Rot}(\mathrm{Aut}(U^+))$. So, since $1_{V(U^+)} \in \Phi$ and $|\Phi(x_\nu)|=|\Phi(w)|=|\Phi|$, clearly $\Phi$ satisfies properties (a) and (b) of Definition \ref{supercard number of common cards}.

(a) It follows from Theorem \ref{max sat set lemma} that there exists a maximum saturating set $X$ of $B_{vw}(U^+)$ such that $\Phi \subseteq X_{\mathrm{Aut}}$. Let $Y$ be any maximum saturating set of $B_{vw}(U^+)$. Suppose that case (i) holds. Then, since it follows from Corollary \ref{Aut rot ref cor 1a} and (i) that $\lambda^{-1}(x_\nu) \in \Phi(x_\nu)$ for all $\lambda \in \mathrm{Ref}(\mathrm{Aut}(U^+))$, we see that $\lambda^{-1}(x_\nu) \in \Phi(x_\nu)$ for all $\lambda \in  \mathrm{Aut}(U^+)$. By Definition \ref{supercard number of common cards}(b), the vertices $\lambda^{-1}(x_\nu)$ are distinct for each $\lambda \in Y$. So $|Y_{\mathrm{Aut}}| \le |\Phi(x_\nu)|=|\Phi|$, and therefore $X_{\mathrm{Aut}} = \Phi$ and thus $X$ is $U^+$-optimum.

Suppose instead that case (ii) holds. Then $\lambda(x_0) \in \Phi(x_0)$ for all $\lambda \in \mathrm{Aut}(U^+)$, and hence $\lambda(w) \in \Phi(w)$ for all $\lambda \in \mathrm{Aut}(U^+)$ as $d_{U^+}(x_0) =3$. By Definition \ref{supercard number of common cards}(b), the vertices $\lambda(w)$ are distinct for each $\lambda \in Y$. So $|Y_{\mathrm{Aut}}| \le |\Phi(w)|=|\Phi|$, and therefore $X_{\mathrm{Aut}} = \Phi$ and thus $X$ is $U^+$-optimum.

(b) We first show that there exists $\pi \in \mathrm {Ref}(\mathrm {Aut}(U^+))$ such that $\Phi(x_\nu) \cap \Phi\pi(x_\nu) =\emptyset$ and $\Phi(w) \cap \Phi\pi(w) =\emptyset$.

Now, since (i) does not hold, we may choose $\sigma \in \mathrm {Ref}(\mathrm {Aut}(U^+))$ such that $\sigma(x_\nu) \not \in \Phi(x_\nu)$. If $\sigma(x_0) \not\in\Phi(x_0)$ then $\sigma(w) \not \in \Phi(w)$, and we simply set $\pi=\sigma$. So suppose that $\sigma(x_0) \in \Phi(x_0)$ and let $x_\alpha=\sigma(x_0)$. By Lemma \ref{X_U cyclic auto subgroup lemma 4}(d), $\mathrm {Ref}(\mathrm {Aut}(U^+))(x_0) \subseteq \Phi(x_0)$, and hence $d_{U^+}(x_0) \ge 4$ since case (ii) does not hold. So $d_{U^+}(x_\alpha) \ge 4$ as $\sigma \in \mathrm {Aut}(U^+)$. We therefore define $\pi \in \mathbf S_{V(U^+)}$ such that $\pi(w)=x_\alpha^2$, $\pi(\sigma^{-1}(x_\alpha^2))=\sigma(w)$, and $\pi(u)=\sigma(u)$ for all other vertices $u$ of $U^+$. Clearly, $\pi \in \mathrm {Ref}(\mathrm{Aut}(U^+))$. Furthermore, $\Phi(x_\nu) \cap \Phi\pi(x_\nu) =\emptyset$ as $\pi(x_\nu) =\sigma(x_\nu)$, and $\Phi(w) \cap \Phi\pi(w) =\emptyset$ by construction.

By Lemma \ref{supercard sunshine iso cards lemma 1}(b), $\lambda^{-1}(x_\nu)=\lambda(x_\nu)$ for all $\lambda \in \mathrm {Ref}(\mathrm {Aut}(U^+))$. So $\{\lambda^{-1}(x_\nu) \mid \lambda \in \Phi\pi\}=\Phi\pi(x_\nu)$ as $\Phi\pi \subseteq \mathrm{Ref}(\mathrm{Aut}(U^+))$. It is now easy to show that the set of permutations $\Phi \cup \Phi\pi$ satisfies property (b) of Definition \ref{supercard number of common cards}. As $1_{V(U^+)} \in \Phi$, it then follows from Theorem \ref{max sat set lemma} that there exists a maximum saturating set $X$ of $B_{vw}(U^+)$ such that $\Phi \subseteq X_{\mathrm{Rot}}$ and $\Phi\pi \subseteq X_{\mathrm{Ref}}$.

Let $Y$ be any maximum saturating set of $B_{vw}(U^+)$. If $\lambda \in Y_{\mathrm{Rot}}$ then $\lambda^{-1}(x_\nu) \in \Phi(x_\nu)$ and if $\lambda \in Y_{\mathrm{Ref}}$ then $\lambda^{-1}(x_\nu) \in \Phi\pi(x_\nu)$. By Definition \ref{supercard number of common cards}(b), the vertices $\lambda^{-1}(x_\nu)$ are distinct for each $\lambda \in Y$. Therefore $|Y_{\mathrm{Rot}}| \le |\Phi |$ and $|Y_{\mathrm{Ref}}| \le |\Phi\pi|$. Hence $X_{\mathrm{Rot}}=\Phi$, $X_{\mathrm{Ref}}=\Phi\pi$ and thus $X$ is $U^+$-optimum. Clearly, $\Phi  \cup \Phi\pi $ is isomorphic to the dihedral group of order $\f{2c}{\delta}$, as $\Phi$ is isormorphic to the cyclic group of order $\f{c}{\delta}$.
\hfill\eot}
\end{thm}
The following lemma concerning the orbits of $X_{\mathrm {Aut}}$ follows from  Definition \ref{supercard number of common cards}, and the fact that the inverses of $X_{\mathrm {Aut}}$ are also in $X_{\mathrm {Aut}}$, since $X_{\mathrm {Aut}}$ is a group.
\begin{lemma}\label{X_U cyclic auto subgroup lemma 6a}
{\rm Let $X$ be the $U^+$-optimum saturating set from Theorem \ref{X_U cyclic auto subgroup thm 6}.
\begin{enumerate}\lettnum\vspace{-4mm}
\item $X_{\mathrm {Aut}}(x_\nu)=\{\lambda^{-1}(x_\nu) \mid \lambda \in X_{\mathrm {Aut}} \}$ and $|X_{\mathrm {Aut}}(x_\nu)|=|X_{\mathrm {Aut}}|$. In addition, $d_{U^+}(x_i)=2$ and $\tau_{U^+}(x_i)=\tau_{U^+}(x_\nu)$ for each vertex $x_i$ in $X_{\mathrm {Aut}}(x_\nu)$.\vspace{-2mm}
\item $|X_{\mathrm {Aut}}(w)|=|X_{\mathrm {Aut}}|$, and each vertex in $X_{\mathrm {Aut}}(w)$ is a leaf adjacent to a vertex of degree $d_{U^+}(x_0)$.\vspace{-2mm}
\item If $\lambda \in X \setminus X_{\mathrm {Aut}}$ then $\lambda^{-1}(x_\nu) \not\in X_{\mathrm {Aut}}(x_\nu)$ and $\lambda(w) \not\in X_{\mathrm {Aut}}(w)$.\vspace{-2mm}
\end{enumerate}\vspace{-4mm}
\hfill\eot}
\end{lemma}
\begin{lemma}\label{X_U cyclic auto subgroup lemma 7}
{\rm Let $X$ be the $U^+$-optimum saturating set from Theorem \ref{X_U cyclic auto subgroup thm 6} and let $B = \{\lambda(w) \mid \lambda \in X \setminus X_{\mathrm {Aut}} \text{ and } d_{U^+}(\lambda(w))=1\}$. Let $A$ be any subset of $B$ such that $\Phi(a_1) \cap \Phi(a_2)=\emptyset$ for all distinct $a_1$ and $a_2$ in $A$. Then
\begin{enumerate}\lettnum\vspace{-4mm}
\item if $X_{\mathrm{Aut}} =\Phi$ then $|X_{\mathrm{Aut}}|  \le \left\lfloor\f{n+1}{\delta+1+|A|}\right\rfloor$;\vspace{-2.5mm}
\item if $X_{\mathrm{Aut}} =\Phi  \cup \Phi\pi$ then $|X_{\mathrm{Aut}}| \le 2\left\lfloor \f{n+1}{\delta+2+|A|}\right\rfloor$.\vspace{-2.5mm}
\end{enumerate}
When $|A|=0$, equality in either case (a) or case (b) holds if and only if $d_1(U^+) = |X_{\mathrm{Aut}}|$.

\textit{Proof} Suppose there exists $a_1 \in A$. By Lemma \ref{X_U cyclic auto subgroup lemma 6a}(c), $a_1 \not\in X_{\mathrm {Aut}}(w)$. So $a_1 \not\in \Phi(w)$, and thus $\Phi(a_1) \cap \Phi(w)=\emptyset$. Similarly, $\Phi(a_1) \cap \Phi\pi(w)=\emptyset$ when $X_{\mathrm {Aut}}=\Phi \cup \Phi\pi$. Therefore $\Phi(a_1) \cap X_{\mathrm {Aut}}(w)=\emptyset$. So, since $\Phi(a_1) \cap \Phi(a_2)=\emptyset$ for all $a_2 \in A \setminus \{a_1\}$, it follows that $d_1(U^+) \ge |A||\Phi|+|X_{\mathrm {Aut}}(w)|$. Hence $d_1(U^+) \ge |A||\Phi|+|X_{\mathrm {Aut}}|$ by Lemma \ref{X_U cyclic auto subgroup lemma 6a}(b). Together with equation (\ref{X_U cyclic auto subgroup lem 7 eqn 1}), this yields $n+1 \ge |X_{\mathrm{Aut}}|+(|A|+\delta)|\Phi|$. Cases (a) and (b) immediately follow. Clearly, when $|A|=0$, equality holds if and only if $d_1(U^+) = |X_{\mathrm{Aut}}|$.
\hfill\eot}
\end{lemma}
\vspace{-6mm}
\section{\normalsize Upper bounds on $b(U,\, T)$}\label{upper bounds on b(U, T)}\vspace{-4mm}
We now obtain upper bounds on $b(U,\, T)$ for the three possible values of $\tau_{U^+}(x_\nu)$, and the two possible cases for $X_{\mathrm{Aut}}$ from Theorem \ref{X_U cyclic auto subgroup thm 6}. We further show that, for each of the six cases, the maximum value of $b(U,\, T)$ is attained if and only if the structure of $U^+$ is as described below.
\begin{enumerate}
\item[(S$0$a)] $n \equiv 2 \pmod 3$, $d_{U^+}(x_i)=3$ when $i \equiv 0 \pmod 2$, and $d_{U^+}(x_i)=2$ otherwise; in addition, $\nu \equiv 1 \pmod 2$.\vspace{-2.5mm}
\item[(S$0$b)] $n \equiv 6 \pmod 7$, $d_{U^+}(x_i)=4$ when $i \equiv 0 \pmod 4$, $d_{U^+}(x_i)=3$ when $i \equiv 2 \pmod 4$, and $d_{U^+}(x_i)=2$ otherwise; in addition, $\nu \equiv 1 \pmod 4$.\vspace{-2.5mm}
\item[(S$1$a)] $n \equiv 3 \pmod 4$, $d_{U^+}(x_i)=3$ when $i \equiv 0 \pmod 3$, and $d_{U^+}(x_i)=2$ otherwise; in addition, $\nu \equiv 1 \pmod 3$.\vspace{-2.5mm}
\item[(S$1$b)] $n \equiv 4 \pmod 5$, $d_{U^+}(x_i)=4$ when $i \equiv 0 \pmod 3$, and $d_{U^+}(x_i)=2$ otherwise; in addition, $\nu \equiv 1 \pmod 3$.\vspace{-2.5mm}
\item[(S$2$a)] $n \equiv 4 \pmod 5$, $d_{U^+}(x_i)=3$ when $i \equiv 0 \pmod 4$, and $d_{U^+}(x_i)=2$ otherwise; in addition, $\nu \equiv 2 \pmod 4$.\vspace{-2.5mm}
\item[(S$2$b)] $n \equiv 6 \pmod 7$, $d_{U^+}(x_i)=4$ when $i \equiv 0 \pmod 5$, and $d_{U^+}(x_i)=2$ otherwise; in addition, either $\nu \equiv 2 \pmod 5$ or $\nu \equiv 3 \pmod 5$.\vspace{-2.5mm}
\end{enumerate}
We note that, $\tau_{U^+}(x_\nu)=0$ when $U^+$ has structure (S$0$a) or (S$0$b); $\tau_{U^+}(x_\nu)=1$ when $U^+$ has structure (S$1$a) or (S$1$b); $\tau_{U^+}(x_\nu)=2$ when $U^+$ has structure (S$2$a) or (S$2$b). We will see that for structures (S$0$a), (S$1$a) and (S$2$a), the corresponding $X_{\mathrm{Aut}}$ is $\Phi$, and for structures (S$0$b), (S$1$b) and (S$2$b), the corresponding $X_{\mathrm{Aut}}$ is $\Phi \cup \Phi\pi$.

Examples of the structures of the six possible supercards are shown in Figures $1$, $2$ and $3$.
\begin{figure}[b]
\centering
\begin{tikzpicture}[thick,scale=1.5]
    \draw \foreach \x in {0,30,...,360}
    {(\x:1) node {} -- (\x+30:1)};
    \draw \foreach \y in {30,90,...,330}
    {(\y:1.5) node {} -- (\y:1.05)};
    \coordinate [label=left:$w$] (y1) at (25:1.9);
    \coordinate [label=left:$x_0$] (y3) at (25:0.9);
	\coordinate [label=left:$x_\nu$] (y2) at (180:0.55);
 \end{tikzpicture}
\hspace{30mm} 
\begin{tikzpicture}[thick,scale=1.5]
    \draw \foreach \x in {0,30,...,360}
    {(\x:1) node {} -- (\x+30:1)};
    \draw \foreach \y in {90,210,...,330}
    {(\y:1.5) node {} -- (\y:1.05)};
    \draw \foreach \y in {15,135,...,255}
    {(\y:1.5) node {} -- (\y+15:1.05)};
    \draw \foreach \y in {45,165,...,285}
    {(\y:1.5) node {} -- (\y-15:1.05)};
    \coordinate [label=left:$w$] (y1) at (12:1.9);
    \coordinate [label=left:$x_0$] (y3) at (25:0.9);
	\coordinate [label=left:$x_\nu$] (y2) at (260:0.75);
 \end{tikzpicture}
\caption{The supercard $U^+$ with structure (S$0$a) when $n=17$ and structure (S$0$b) when $n=20$.} \label{Fig Thm 6a}
\end{figure}
\begin{figure}
\centering
\begin{tikzpicture}[thick,scale=1.5]
    \draw \foreach \x in {0,30,...,360}
    {(\x:1) node {} -- (\x+30:1)};
    \draw \foreach \y in {30,120,...,300}
    {(\y:1.5) node {} -- (\y:1.05)};
    \coordinate [label=left:$w$] (y1) at (25:1.9);
    \coordinate [label=left:$x_0$] (y3) at (25:0.9);
	\coordinate [label=left:$x_\nu$] (y2) at (180:0.55);
 \end{tikzpicture}
\hspace{30mm} 
\begin{tikzpicture}[thick,scale=1.5]
    \draw \foreach \x in {0,30,...,360}
    {(\x:1) node {} -- (\x+30:1)};
    \draw \foreach \y in {15,105,...,315}
    {(\y:1.5) node {} -- (\y+15:1.05)};
    \draw \foreach \y in {45,135,...,345}
    {(\y:1.5) node {} -- (\y-15:1.05)};
    \coordinate [label=left:$w$] (y1) at (12:1.9);
    \coordinate [label=left:$x_0$] (y3) at (25:0.9);
	\coordinate [label=left:$x_\nu$] (y2) at (180:0.55);
 \end{tikzpicture}
\caption{The supercard $U^+$ with structure (S$1$a) when $n=15$ and structure (S$1$b) when $n=19$.}\label{Fig Thm 6b}
\end{figure}
\begin{figure}
\centering
\begin{tikzpicture}[thick,scale=1.5]
    \draw \foreach \x in {0,30,...,360}
    {(\x:1) node {} -- (\x+30:1)};
    \draw \foreach \y in {30,150,...,270}
    {(\y:1.5) node {} -- (\y:1.05)};
    \coordinate [label=left:$w$] (y1) at (25:1.9);
    \coordinate [label=left:$x_0$] (y3) at (25:0.9);
	\coordinate [label=left:$x_\nu$] (y2) at (220:0.6);
 \end{tikzpicture}
\hspace{30mm} 
\begin{tikzpicture}[thick,scale=1.5]
    \draw \foreach \x in {0,24,...,360}
    {(\x:1) node {} -- (\x+24:1)};
    \draw \foreach \y in {14,134, 254}
    {(\y:1.5) node {} -- (\y+10:1.05)};
    \draw \foreach \y in {34,154, 274}
    {(\y:1.5) node {} -- (\y-10:1.05)};
    \coordinate [label=left:$w$] (y1) at (12:1.9);
    \coordinate [label=left:$x_0$] (y3) at (20:0.9);
	\coordinate [label=left:$x_\nu$] (y2) at (230:0.60);
 \end{tikzpicture}
\caption{The supercard $U^+$ with structure (S$2$a) when $n=14$ and structure (S$2$b) when $n=20$.}\label{Fig Thm 6c}
\end{figure}

We first prove that if $U^+$ has any of these structures then $B_{vw}^U(U^+)=\mathrm{Aut}(U^+)$, and specify the value of $b(U,\, T)$ in each case.
\begin{lemma}\label{sunshine max cards auto lemma 1}
{\rm Suppose that $U^+$ has any of the six structures (S$0$a) through to (S$2$b). Then then there exists a non-trivial rotation in $\mathrm{Aut}(U^+)$ and $B_{vw}^U(U^+)=\mathrm{Aut}(U^+)$.

\textit{Proof} It is easy to see that if $U^+$ has any of the six structures, then there exists a non-trivial rotation in $\mathrm{Aut}(U^+)$. So $B^U_{vw}(U^+)=\mathrm{Aut}(U^+)$ by Lemma \ref{sunshine supercard one rot auto all rot autos lemma}, in each case.
\hfill\eot}
\end{lemma}
\vspace{-8mm}
\newpage
\begin{lemma}\label{sunshine max cards auto lemma 1a}
{\rm Let $n \ge 14$. Suppose that $U^+$ has any of the six structures (S$0$a) through to (S$2$b).
\begin{enumerate}\lettnum\vspace{-4mm}
\item If $U^+$ has structure (S$0$a) then $b(U,\, T) = \f{n+1}{3}+2$.\vspace{-2.5mm}
\item If $U^+$ has structure (S$0$b) then $b(U,\, T) = \f{2(n+1)}{7}+1$.\vspace{-2.5mm}
\item If $U^+$ has structure (S$1$a) then $b(U,\, T) = \f{n+1}{4}+1$.\vspace{-2.5mm}
\item If $U^+$ has structure (S$1$b) then $b(U,\, T) = \f{2(n+1)}{5}$.\vspace{-2.5mm}
\item If $U^+$ has structure (S$2$a) then $b(U,\, T) = \f{n+1}{5}$.\vspace{-2.5mm}
\item If $U^+$ has structure (S$2$b) then $b(U,\, T) = \f{2(n+1)}{7}$.
\end{enumerate}
\textit{Proof} In each of the six cases, there exists a non-trivial rotation in $\mathrm{Aut}(U^+)$ by Lemma \ref{sunshine max cards auto lemma 1}, so we therefore apply all the results from Section \ref{B(U^+) = U^+}. In particular, we may define $\Phi$ and $\delta$ as in Definition \ref{Phi and delta}. 

Suppose that $U^+$ has structure (S$0$a), (S$1$a) or (S$2$a). Then clearly, case (ii) of Theorem \ref{X_U cyclic auto subgroup thm 6} holds. So there exists a $U^+$-optimum saturating set $X$ such that $X_{\mathrm {Aut}}=\Phi$.

Suppose instead that $U^+$ has structure (S$0$b), (S$1$b) or (S$2$b). Then, it is easy to see that there exists $\lambda \in \mathrm{Ref}(\mathrm{Aut}(U^+))$ such that $\lambda(x_0)=x_0$, and to verify that $\lambda(x_\nu)$, i.e., $x_{c-\nu}$, is not in $\Phi(x_\nu)$. So case (i) of Theorem \ref{X_U cyclic auto subgroup thm 6} cannot hold. Moreover, case (ii) of Theorem \ref{X_U cyclic auto subgroup thm 6} cannot hold as $d_{U^+}(x_0) \ge 4$. Therefore, there exists $\pi \in \mathrm{Ref}(\mathrm{Aut}(U^+))$ and a $U^+$-optimum saturating set $X$ such that $X_{\mathrm {Aut}}=\Phi \cup \Phi\pi$.

(a) Since $d_1(U^+)=|\Phi|$, it follows from equation (\ref{X_U cyclic auto subgroup lem 7 eqn 1}) that $|X_{\mathrm{Aut}}|= |\Phi|=\f{n+1}{3}$. It is easy to see that $U-x_2 \cong U-x_{c-2} \cong T-x_{\nu+1} \cong  T-x_{\nu-1}$, and to check using Lemma \ref{supercard sunshine lemma 1a} that $\lambda(w) \in \{x_{\nu-1},\,  x_{\nu+1}\}$ for any $\lambda \in X \setminus X_{\mathrm{Aut}}$. Therefore $|X \setminus X_{\mathrm{Aut}}|=2$, and $b(U,\, T) = |X|=\f{n+1}{3}+2$.

(b) Since $d_1(U^+)=3|\Phi|$, it follows from equation (\ref{X_U cyclic auto subgroup lem 7 eqn 1}) that $|X_{\mathrm{Aut}}|= 2|\Phi|=\f{2(n+1)}{7}$. It is easy to see that $U-x_2 \cong U-x_{c-2} \cong T-x_{\nu+1}$, and to check using Lemma \ref{supercard sunshine lemma 1a} that $\lambda(w) = x_{\nu+1}$ for any $\lambda \in X \setminus X_{\mathrm{Aut}}$. Therefore $|X \setminus X_{\mathrm{Aut}}|=1$, and $b(U,\, T) = |X|=\f{2(n+1)}{7}+1$.

(c) Since $d_1(U^+)=|\Phi|$, it follows from equation (\ref{X_U cyclic auto subgroup lem 7 eqn 1}) that $|X_{\mathrm{Aut}}|= |\Phi|=\f{n+1}{4}$. It is easy to see that $U-x_2 \cong U-x_{c-2} \cong T-x_{\nu+1}$, and to check that to check using Lemma \ref{supercard sunshine lemma 1a} that $\lambda(w) = x_{\nu+1}$ for any $\lambda \in X \setminus X_{\mathrm{Aut}}$. Therefore $|X \setminus X_{\mathrm{Aut}}|=1$, and $b(U,\, T) = |X|=\f{n+1}{4}+1$.

For (d), (e), and (f), it is easy to check using Lemma \ref{supercard sunshine lemma 1a} that $\lambda(w) \in V(\mathrm{skel}(U^+))$ for any $\lambda \in X \setminus X_{\mathrm{Aut}}$. Furthermore, it is easy to see that $U-x_i \not\cong T-x_j$, for all $i,\, j$, and hence $X = X_{\mathrm{Aut}}$ in each case.

(d) Since $d_1(U^+)=2|\Phi|$, it follows from equation (\ref{X_U cyclic auto subgroup lem 7 eqn 1}) that $|X_{\mathrm{Aut}}|= 2|\Phi|=\f{2(n+1)}{5}$. Therefore $b(U,\, T) = |X|=\f{2(n+1)}{5}$.

(e) Since $d_1(U^+)=|\Phi|$, it follows from equation (\ref{X_U cyclic auto subgroup lem 7 eqn 1}) that $|X_{\mathrm{Aut}}|= |\Phi|=\f{n+1}{5}$. Therefore $b(U,\, T) = |X|=\f{n+1}{5}$.

(f) Since $d_1(U^+)=2|\Phi|$, it follows from equation (\ref{X_U cyclic auto subgroup lem 7 eqn 1}) that $|X_{\mathrm{Aut}}|=2|\Phi|=\f{2(n+1)}{7}$. Therefore $b(U,\, T) = |X|=\f{2(n+1)}{7}$.
\hfill\eot}
\end{lemma}
For each of the following four lemmas, we assume that {\em there exists a non-trivial rotation in $\mathrm{Aut}(U^+)$, and we let $X$ be the $U^+$-optimum saturating set of $B_{vw}(U^+)$ from Theorem \ref{X_U cyclic auto subgroup thm 6}}. We recall that $\widetilde{X}$ is the subset of $X$ containing those permuations $\lambda$ such that $\lambda(w)$ is a leaf of $U^+$ and a $d$-leaf of $T$. We note that $X_U=X_{\mathrm{Aut}}$ by Lemma \ref{sunshine supercard one rot auto all rot autos lemma}.
\begin{lemma}\label{sunshine max cards auto lemma 2}
{\rm Let $n \ge 56$ and suppose that $\tau_{U^+}(x_\nu)=0$.
\begin{enumerate}\lettnum\vspace{-4mm}
\item If $X_{\mathrm{Aut}}=\Phi$ then $b(U,\, T) =|X| \le \f{n+1}{3}+2$. Furthermore, equality holds if and only if $U^+$ has structure (S$0$a).\vspace{-2.5mm}
\item If $X_{\mathrm{Aut}}=\Phi \cup \Phi\pi$ then $b(U,\, T) =|X| \le  \f{2(n+1)}{7}+1$. Furthermore, equality holds if and only if $U^+$ has structure (S$0$b).\vspace{-2.5mm}
\end{enumerate}
\textit{Proof} It follows from Corollary \ref{supercard sunshine cor 7a} that $\widetilde{X} \subseteq X_{\mathrm{Aut}}$, so $X \setminus X_{\mathrm{Aut}} \subseteq  X \setminus \widetilde{X}$. Therefore $|X \setminus X_{\mathrm{Aut}}| \le 4$ by Corollary \ref{supercard sunshine cor 4c}(a).

(a) By Lemma \ref{X_U cyclic auto subgroup lemma 7}(a), $|X_{\mathrm{Aut}}| \le \left\lfloor \f{n+1}{\delta+1} \right\rfloor$. Simple calculations then show that the bound holds for $\delta \ge 3$ with strict inequality. 

So suppose that $\delta =2$. If $d_{U^+}(x_0) \ge 4$, then $d_1(U^+) \ge 2|\Phi|$, so $|X_{\mathrm{Aut}}| =|\Phi| \le \left\lfloor \f{n+1}{4} \right\rfloor$ by (\ref{X_U cyclic auto subgroup lem 7 eqn 1}). Simple calculations then show that the bound holds with strict inequality in this case. On the other hand, if $d_{U^+}(x_0) =3$ then $U^+$ has structure (S$0$a), so $b(U,\, T)=\f{n+1}{3}+2$ by Lemma \ref{sunshine max cards auto lemma 1a}(a). Therefore, the bound holds in all cases, and is attained if and only if $U^+$ has structure (S$0$a). 

(b) By Lemma \ref{X_U cyclic auto subgroup lemma 6a}(a), $U^+$ contains $2|\Phi|$ vertices $x_i$ with $\tau_{U^+}(x_i)=2$. As $\tau_{U^+}(x_\nu)=0$, this implies that $\delta \ge 4$. In addition, it is easy to see by inspection, that $d_1(U^+) \ge 3|\Phi|$ when $\delta \ge 5$.

Suppose that $\delta \ge 5$. Then $|X_{\mathrm{Aut}}| \le 2\left\lfloor \f{n+1}{\delta+3} \right\rfloor$ by (\ref{X_U cyclic auto subgroup lem 7 eqn 1}). Simple calculations then show that the bound holds with strict inequality when $\delta \ge 6$ or $|X \setminus X_{\mathrm{Aut}}| = 3$. So suppose that $\delta=5$ and $|X \setminus X_{\mathrm{Aut}}| = 4$. In this case, it is easy to see from Lemma \ref{supercard sunshine lemma 4b}(a) that $\{x^1_{\nu-1}, x^1_{\nu+1}\}$ is the set $B$ of leaves defined in Lemma \ref{X_U cyclic auto subgroup lemma 7}. Moreover, since $\delta \ge 5$, we may clearly put $A=B$ in the lemma. Hence $|X_{\mathrm{Aut}}| \le 2\left\lfloor \f{n+1}{9} \right\rfloor$ by Lemma \ref{X_U cyclic auto subgroup lemma 7}(b), and again simple calculations show that the bound holds with strict inequality.

So suppose that $\delta = 4$. Then $d_{U^+}(x_{\nu-1}) > d_{U^+}(x_{\nu+1}) \ge 3$ and $d_{U^+}(x_{\nu+2})=2$. Now, if $d_{U^+}(x_{\nu+1}) \ge 4$, then $d_1(U^+) \ge 5|\Phi|$ and $|X_{\mathrm{Aut}}| \le 2\left\lfloor \f{n+1}{9} \right\rfloor$ by (\ref{X_U cyclic auto subgroup lem 7 eqn 1}).  Simple calculations then show that the bound holds with strict inequality in this case. On the other hand, if $d_{U^+}(x_{\nu-1}) =4$ and $d_{U^+}(x_{\nu-1}) =3$, then $U^+$ has structure (S$0$b), so $b(U,\, T)=\f{2(n+1)}{7}+1$ by Lemma \ref{sunshine max cards auto lemma 1a}(b). Hence the bound holds in all cases, and is attained if and only if $U^+$ has structure (S$0$b). 
\hfill\eot}
\end{lemma} 
\begin{lemma}\label{sunshine max cards auto lemma 3}
{\rm Let $n \ge 60$ and suppose that $\tau_{U^+}(x_\nu)=2$.
\begin{enumerate}\lettnum\vspace{-4mm}
\item If $X_{\mathrm{Aut}}=\Phi$ then $b(U,\, T) \le \f{n+1}{5} $. Furthermore, equality holds in the bound if and only if $U^+$ has structure (S$2$a).\vspace{-2.5mm}
\item If $X_{\mathrm{Aut}}=\Phi \cup \Phi\pi$ then $b(U,\, T) \le \f{2(n+1)}{7}$. Furthermore, equality holds in the bound if and only if $U^+$ has structure (S$2$b).\vspace{-2.5mm}
\end{enumerate}
\textit{Proof} It follows from Corollary \ref{supercard sunshine cor 7a} that $\widetilde{X} \subseteq X_{\mathrm{Aut}}$, so $X \setminus X_{\mathrm{Aut}} \subseteq  X \setminus \widetilde{X}$. Therefore $|X \setminus X_{\mathrm{Aut}}| \le 2$ by Corollary \ref{supercard sunshine cor 4c}(c).

(a) By Lemma \ref{X_U cyclic auto subgroup lemma 7}(a), $|X_{\mathrm{Aut}}| \le \left\lfloor \f{n+1}{\delta+1} \right\rfloor$. Simple calculations then show that the bound holds for $\delta \ge 5$ with strict inequality. 

So suppose that $\delta \le 4$. By Lemma \ref{X_U cyclic auto subgroup lemma 6a}(a), $U^+$ contains $|\Phi|$ vertices $x_i$ with $\tau_{U^+}(x_i)=2$. This implies that $\delta=4$, and that every cut-vertex of $U^+$ is in $\Phi(x_0)$. Now, if $d_{U^+}(x_0) \ge 4$ then $d_1(U^+) \ge 2|\Phi|$, so $|X_{\mathrm{Aut}}| \le \left\lfloor \f{n+1}{6} \right\rfloor$ by (\ref{X_U cyclic auto subgroup lem 7 eqn 1}). Simple calculations then show that the bound holds with strict inequality in this case. On the other hand, if $d_{U^+}(x_0) =3$ then $U^+$ has structure (S$2$a), so $b(U,\, T)=\f{n+1}{5}$ by Lemma \ref{sunshine max cards auto lemma 1a}(e). Therefore, the bound holds in all cases, and is attained if and only if $U^+$ has structure (S$2$a).

(b) By Lemma \ref{X_U cyclic auto subgroup lemma 7}(b), $|X_{\mathrm{Aut}}| \le 2\left\lfloor \f{n+1}{\delta+2} \right\rfloor$.  Simple calculations then show that the bound holds for $\delta \ge 6$ with strict inequality.

So suppose that $\delta \le 5$. By Lemma \ref{X_U cyclic auto subgroup lemma 6a}(a), $U^+$ contains $2|\Phi|$ vertices $x_i$ with $\tau_{U^+}(x_i)=2$. This implies that $\delta=5$, and that every cut-vertex of $U^+$ is in $\Phi(x_0)$. Since  $U^+$ contains at least $2|\Phi|$ leaves by Lemma \ref{X_U cyclic auto subgroup lemma 6a}(b), clearly $d_{U^+}(x_0) \ge 4$. Now, if $d_{U^+}(x_0) \ge 5$ then $d_1(U^+) \ge 3|\Phi|$, so $|X_{\mathrm{Aut}}| \le 2\left\lfloor \f{n+1}{8} \right\rfloor$ by (\ref{X_U cyclic auto subgroup lem 7 eqn 1}). Simple calculations then show that the bound holds with strict inequality in this case. On the other hand, if $d_{U^+}(x_0) =4$ then $U^+$ has structure (S$2$b), so $b(U,\, T)=\f{2(n+1)}{7}$ by Lemma \ref{sunshine max cards auto lemma 1a}(f). Therefore, the bound holds in all cases, and is attained if and only if $U^+$ has structure (S$2$b).
\hfill\eot}
\end{lemma} 
Since $\widetilde{X}$ may not be contained in $X_{\mathrm{Aut}}$ when $\tau_{U^+}(x_\nu)=1$, we need an auxillary result in this case.
\begin{lemma}\label{sunshine max cards auto lemma 4}
{\rm Suppose that $\lambda \in \widetilde{X} \setminus X_{\mathrm{Aut}}$ and that $\delta \le \f{c-3}{2}$. Then $\lambda(w)$ is adjacent to $x_{\nu+2}$ and $\lambda^{-1}(x_\nu) \in \{x_2,\, x_{c-2}\}$.

\textit{Proof} Let $x_\mu=\lambda^{-1}(x_\nu)$. By Corollary \ref{supercard sunshine cor 7a},  $\tau_{U^+}(x_\nu)=1$ and $\{\lambda(x_{\mu+2}),\, \lambda(x_{\mu-2})\} = \{x_{\nu+3},\, x_{\nu-1}\}$ as $\lambda \in \widetilde{X} \setminus X_{\mathrm{Aut}}$. Moreover, by Lemma \ref{supercard sunshine lemma 5a_2}, either
\begin{enumerate}\lettnum\vspace{-4mm}
\item $\mathrm{skel}(U-x_\mu)$ is $x_{\mu+1}x_{\mu+2} \ldots x_{\mu-2}$ and $\lambda(x_i)= x_{(\nu-\mu+1)+i}$ for all $x_i$ in $V(\mathrm{skel}(U-x_\mu))$, or \vspace{-2.5mm}
\item $\mathrm{skel}(U-x_\mu)$ is $x_{\mu+2}x_{\mu+3} \ldots x_{\mu-1}$ and $\lambda(x_i)= x_{(\nu+\mu+1)-i}$ for all $x_i$ in $V(\mathrm{skel}(U-x_\mu))$. 
\end{enumerate} \vspace{-2.5mm}
We recall that $d_{U^+}(x_i)=d_{U^+}(x_{i+\delta})=d_{U^+}(x_{i-\delta})$ for all $x_i$, as $\phi(x_i)=x_{i+\delta}$.

Case (a): We first note that $\{x_{\mu+1+\delta}, \, x_{\mu+1+2\delta}, \, x_{\mu-2-\delta}\} \subseteq V(\mathrm{skel}(U-x_\mu))$ as $\delta \le \f{c-3}{2}$. Suppose that $\lambda(w)$ is not adjacent to $x_{\nu+2}$. Then
\[
d_{U^+}(x_{\mu+1}) > d_{U-x_\mu}(x_{\mu+1})=d_{T-\lambda(w)}(\lambda(x_{\mu+1}))=d_{T-\lambda(w)}(x_{\nu+2})=d_{U^+}(x_{\nu+2}).
\]
So $d_{U^+}(x_{\mu+1+\delta}) > d_{U^+}(x_{\nu+2+\delta})$ and $d_{U^+}(x_{\mu+1+2\delta}) > d_{U^+}(x_{\nu+2+2\delta})$. However, since $\lambda(x_{\mu+1+\delta}) = x_{\nu+2+\delta}$ and $\lambda(x_{\mu+1+2\delta})= x_{\nu+2+2\delta}$, we see that $d_{U-x_\mu}(x_{\mu+1+\delta}) = d_{T-\lambda(w)}(x_{\nu+2+\delta})$ and \\$d_{U-x_\mu}(x_{\mu+1+2\delta}) = d_{T-\lambda(w)}(x_{\nu+2+2\delta})$. Since $x_\mu$ is not adjacent to $x_{\mu+1+\delta}$ or $x_{\mu+1+2\delta}$, it follows that $w$ must be adjacent to both of these vertices. This is impossible since $w$ is a leaf. Therefore $\lambda(w)$ is adjacent to $x_{\nu+2}$.

Suppose now that $x_\mu$ is not $x_2$. Then $w$ is not adjacent to $x_{\mu-2}$, and thus 
\[
d_{U^+}(x_{\nu-1}) > d_{T-\lambda(w)}(x_{\nu-1})=d_{U-x_\mu}(\lambda^{-1}(x_{\nu-1}))=d_{U-x_\mu}(x_{\mu-2})=d_{U^+}(x_{\mu-2}).
\]
So $d_{U^+}(x_{\nu-1-\delta}) > d_{U^+}(x_{\mu-2-\delta})$. However, since $\lambda(x_{\mu-2-\delta}) = x_{\nu-1-\delta}$, it follows that \\$d_{U-x_\mu}(x_{\nu-1-\delta}) = d_{T-\lambda(w)}(x_{\mu-2-\delta})$, which is impossible as neither $\lambda(w)$ nor $x_\nu$ are adjacent to $x_{\nu-1-\delta}$. Therefore $x_\mu$ is $x_2$.

Case (b) can be proved in a similar manner by replacing $x_{\mu+k}$ by $x_{\mu-k}$ for each $k$, and vice versa, and also substituting $x_{c-2}$ for $x_2$.
\hfill\eot}
\end{lemma}  
\begin{lemma}\label{sunshine max cards auto lemma 3}
{\rm Let $n \ge 48$ and suppose that $\tau_{U^+}(x_\nu)=1$.
\begin{enumerate}\lettnum\vspace{-4mm}
\item If $X_{\mathrm{Aut}}=\Phi$ then $b(U,\, T) \le \f{n+1}{4}+1$. Furthermore, equality holds if and only if $U^+$ has structure (S$1$a).\vspace{-2.5mm}
\item If $X_{\mathrm{Aut}}=\Phi \cup \Phi\pi$ then $b(U,\, T) \le \f{2(n+1)}{5}$. Furthermore, equality holds if and only if and $U^+$ has structure (S$1$b).\vspace{-2.5mm}
\end{enumerate}
\textit{Proof} We first note that, $\delta \ge 3$, as $\tau_{U^+}(x_\nu) \ge 1$. When $\delta=3$, it is easy to see that there exists a $\psi \in \mathrm {Ref}(\mathrm{Aut}(U^+))$ such that $\psi(x_\nu)=x_{\nu+1} \not\in \Phi(x_\nu)$, and $\psi(x_0) \in \Phi(x_0)$. Thus $\mathrm {Ref}(\mathrm{Aut}(U^+))(x_\nu) \not \subseteq \Phi(x_\nu)$, whereas $\mathrm {Ref}(\mathrm{Aut}(U^+))(x_0) \subseteq \Phi(x_\nu)$ by Lemma \ref{X_U cyclic auto subgroup lemma 4}(e).

Suppose that $\delta \ge \f{c-2}{2}$, so $|\Phi| =\f{c}{\delta} \le 2$. Then, since $|X_{\mathrm{Aut}}|=\chi(U^+)$, it follows from Theorem \ref{supercard sunshine thm 7e} that $b(U,\, T) \le 9$ when $X_{\mathrm{Aut}} =\Phi$, and $b(U,\, T) \le 13$ when $X_{\mathrm{Aut}} =\Phi \cup \Phi\pi$. Simple calculations show that both bounds hold with strict inequality.

We therefore assume that $\delta \le \f{c-3}{2}$. Now, if $\lambda \in \widetilde{X} \setminus X_{\mathrm{Aut}}$ then $\lambda^{-1}(x_\nu) \in \{x_2,\, x_{c-2}\}$ by Lemma \ref{sunshine max cards auto lemma 4}; hence $|\widetilde{X} \setminus X_{\mathrm{Aut}}| \le 2$ by Definition \ref{supercard number of common cards}(b). Therefore, since $|X \setminus \widetilde{X}| \le 4$ by Corollary \ref{supercard sunshine cor 4c}(b), clearly $|X \setminus X_{\mathrm{Aut}}| \le 6$.

Let $B$ be the set of leaves defined in Lemma \ref{X_U cyclic auto subgroup lemma 7}. Then, by Lemma \ref{supercard sunshine lemma 4b}(b) and Lemma \ref{sunshine max cards auto lemma 4}, $B \subseteq \{x^1_{\nu-1}, \, x_{\nu+2}^j,\, \,x_{\nu+2}^k\}$, for some $j$, $k$. On the other hand, if $\lambda \in X \setminus X_{\mathrm{Aut}}$ but $\lambda(w) \not \in B$ then $\lambda(w) \in \{x_{\nu-1},\, x_{\nu+1}, x_{\nu+2}\}$. It follows that if $|X \setminus X_{\mathrm{Aut}}| \ge 4$, then there exists a subset $A$ of $B$ with $|A| \ge 1 $ that satisfies the conditions of Lemma \ref{X_U cyclic auto subgroup lemma 7}. In addition, if $\delta \ge 4$ and $|X \setminus X_{\mathrm{Aut}}| \ge 5$, it is easy to see that there is some subset $A$ of $B$ with $|A| \ge 2$ that satisfies the conditions of Lemma \ref{X_U cyclic auto subgroup lemma 7}.

(a) Suppose that $\delta \ge 4$. Then, by Lemma \ref{X_U cyclic auto subgroup lemma 7}(a), $|X_{\mathrm{Aut}}| \le \left\lfloor \f{n+1}{5} \right\rfloor$ when $|X \setminus X_{\mathrm{Aut}}| \le 3$, $|X_{\mathrm{Aut}}| \le \left\lfloor \f{n+1}{6} \right\rfloor$  when $|X \setminus X_{\mathrm{Aut}}| =4$, and $|X_{\mathrm{Aut}}| \le \left\lfloor \f{n+1}{7} \right\rfloor$ otherwise. Simple calculations then show that the bound holds with strict inequality.

So, suppose that $\delta =3$. Then, since $X_{\mathrm{Aut}} = \Phi$ and $\mathrm {Ref}(\mathrm{Aut}(U^+))(x_\nu) \not \subseteq \Phi(x_\nu)$, it follows from Theorem \ref{X_U cyclic auto subgroup thm 6} that $d_{U^+}(x_0)=3$. So $U^+$ has structure (S$1$a) and $b(U,\, T)=\f{n+1}{4}+1$ by Lemma \ref{sunshine max cards auto lemma 1a}(c). Therefore, the bound holds in all cases, and is attained if and only if $U^+$ has structure (S$1$a). 

(b) Suppose that $\delta \ge 4$. Then, by Lemma \ref{X_U cyclic auto subgroup lemma 7}(b), $|X_{\mathrm{Aut}}| \le 2\left\lfloor \f{n+1}{6} \right\rfloor$ when $|X \setminus X_{\mathrm{Aut}}| \le 3$, $|X_{\mathrm{Aut}}| \le 2\left\lfloor \f{n+1}{7} \right\rfloor$  when $|X \setminus X_{\mathrm{Aut}}| =4$, and $|X_{\mathrm{Aut}}| \le 2\left\lfloor \f{n+1}{8} \right\rfloor$ otherwise. Simple calculations then show that the bound holds with strict inequality.

So, suppose that $\delta =3$. Then, since $X_{\mathrm{Aut}} =  \Phi \cup \Phi\pi$, $\mathrm {Ref}(\mathrm{Aut}(U^+))(x_\nu) \not \subseteq \Phi(x_\nu)$ and $\mathrm {Ref}(\mathrm{Aut}(U^+))(x_0) \subseteq \Phi(x_0)$, it follows from Theorem \ref{X_U cyclic auto subgroup thm 6} that $d_{U^+}(x_0) \ge 4$. Now,  by Lemma \ref{X_U cyclic auto subgroup lemma 6a}(a), $U^+$ contains $|X_{\mathrm{Aut}}(x_\nu)|$ vertices $x_i$ with $d_{U^+}(x_i)=2$. As $d_{U^+}(x_0) \ge 4$, it then follows from Lemma \ref{X_U cyclic auto subgroup lemma 6a}(c) that if $\lambda \in X \setminus X_{\mathrm{Aut}}$ then $\lambda^{-1}(x_\nu)$ must be a cut-vertex of $U$. It is easy to see by inspection that there can be no such $\lambda$, and therefore $X = X_{\mathrm{Aut}}$. The bound then holds by Lemma \ref{X_U cyclic auto subgroup lemma 7}(b) with $|A|=0$, with equality holding if and only if $d_1(U^+)=2|\Phi|$, i.e., when $d_{U^+}(x_0) = 4$. In this case, $U^+$ has structure ($1$b), and $b(U,\, T)=\f{2(n+1)}{5}$ by Lemma \ref{sunshine max cards auto lemma 1a}(d). Therefore, the bound holds in all cases, and is attained if and only if $U^+$ has structure (S$1$b). 

\vspace{-2.5mm}
\hfill\eot}
\end{lemma} 
By combining the above results with Theorem \ref{X_U not autos thm final}, we finally obtain a bound on $b(U,\, T)$ which holds in all cases.
\begin{thm}\label{final suns and cats thm 1}
{\rm Let $U$ be a sunshine graph and $T$ be a caterpillar of order $n$, where $n \ge 60$. Suppose there exists a sunshine graph $U^+$ that is a supercard of $U$ and $T$ such that $\mathrm{Aut}(U^+)$ contains a non-trivial rotation. Then $b(U,\, T) \le \f{2(n+1)}{5}$, with equality if and only if $U^+$ has structure (S$1$b), in which case $n \equiv 4 \pmod 5$. Moreover, in all other cases, $b(U,\, T) \le \f{n+1}{3}+2$.

\textit{Proof} This follows immediately from Lemmas \ref{sunshine max cards auto lemma 1} to \ref{sunshine max cards auto lemma 3}.
\hfill\eot}
\end{thm}
We note that, with more work, we can show that this bound holds for smaller values of $n$ (this is relatively straightforward for $n \ge 35$). However, since the proofs are slightly technical, in the interests of brevity, we have not included them in this paper.
\begin{thm}\label{final suns and cats thm 2}
{\rm Let $U$ be a sunshine graph and $T$ be a caterpillar, where $n \ge 62$. Then $b(U,\, T) \le \f{2(n+1)}{5}$, with equality if and only if there is a supercard of $U$ and $T$ that has structure (S$1$b), in which case $n \equiv 4 \pmod 5$. Moreover, in all other cases, $b(U,\, T) \le \f{n+1}{3}+2$ when $n \ge 74$.

\textit{Proof} We may clearly assume that $b(U,\, T) \ge 10$. By Lemma \ref{sun cat lemma new 2},  there exists a supercard $U^+$ of $U$ and $T$ that is a sunshine graph, and we may choose $U^+$  to have the largest possible value of $\chi(U^+)$ over all supercards of $U$ and $T$ that are sunshine graphs. If $B_U(U^+) = \mathrm{Aut}(U^+)$  then $\mathrm{Aut}(U^+)$ contains a non-trivial rotation by Corollary \ref{sunshine supercard two distinct auto refs all rot autos cor 1}(c), so the results follow from Theorem \ref{final suns and cats thm 1}. When $B_U(U^+) \ne \mathrm{Aut}(U^+)$, they follow from Theorem \ref{X_U not autos thm final} by straightforward calculation.
\hfill\eot}
\end{thm}

\end{document}